# Where Did Combinators Come From? Hunting the Story of Moses Schönfinkel

Stephen Wolfram*

*Combinators were a key idea in the development of mathematical logic and the emergence of the concept of universal computation. They were introduced on December 7, 1920, by Moses Schönfinkel. This is an exploration of the personal story and intellectual context of Moses Schönfinkel, including extensive new research based on primary sources.*

## December 7, 1920

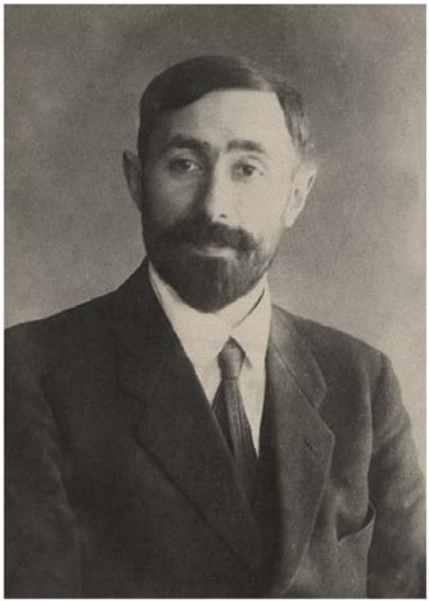
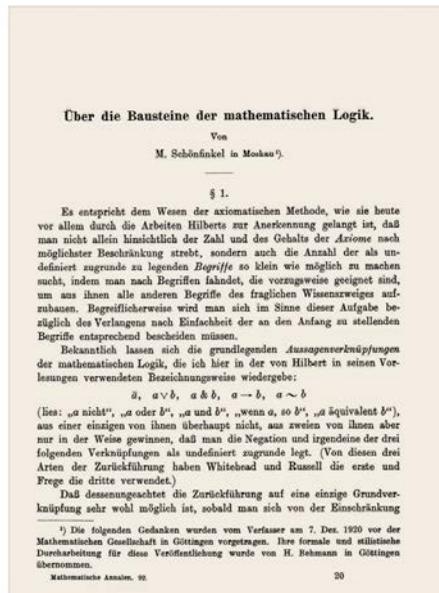

On Tuesday, December 7, 1920, the Göttingen Mathematics Society held its regular weekly meeting—at which a 32-year-old local mathematician named Moses Schönfinkel with no known previous mathematical publications gave a talk entitled "Elemente der Logik" ("Elements of Logic").

A hundred years later what was presented in that talk still seems in many ways alien and futuristic—and for most people almost irreducibly abstract. But we now realize that that talk gave the first complete formalism for what is probably the single most important idea of this past century: the idea of universal computation.

Sixteen years later would come Turing machines (and lambda calculus). But in 1920 Moses Schönfinkel presented what he called "building blocks of logic"—or what we now call "combinators"—and then proceeded to show that by appropriately combining them one could effectively define any function, or, in modern terms, that they could be used to do universal computation.

Looking back a century it's remarkable enough that Moses Schönfinkel conceptualized a formal system that could effectively capture the abstract notion of computation. And it's more remarkable still that he formulated what amounts to the idea of universal computation, and showed that his system achieved it.

But for me the most amazing thing is that not only did he invent the first complete formalism for universal computation, but his formalism is probably in some sense minimal. I've personally spent years trying to work out just how simple the structure of systems that support universal computation can be—and for example with Turing machines it took from 1936 until 2007 for us to find the minimal case.

But back in his 1920 talk Moses Schönfinkel—presenting a formalism for universal computation for the very first time—gave something that is probably already in his context minimal.

Moses Schönfinkel described the result of his 1920 talk in an 11-page paper published in 1924 entitled "Über die Bausteine der mathematischen Logik" ("On the Building Blocks of Mathematical Logic"). The paper is a model of clarity. It starts by saying that in the "axiomatic method" for mathematics it makes sense to try to keep the number of "fundamental notions" as small as possible. It reports that in 1913 Henry Sheffer managed to show that basic logic requires only one connective, that we now call ɴᴀɴᴅ. But then it begins to go further. And already within a couple of paragraphs it's saying that "We are led to [an] idea, which at first glance certainly appears extremely bold". But by the end of the introduction it's reporting, with surprise, the big news: "It seems to me remarkable in the extreme that the goal we have just set can be realized… [and]; as it happens, it can be done by a reduction to three fundamental signs".

Those "three fundamental signs", of which he only really needs two, are what we now call the *S* and *K* combinators (he called them *S* and *C*). In concept they're remarkably simple, but their actual operation is in many ways brain-twistingly complex. But there they were—already a century ago—just as they are today: minimal elements for universal computation, somehow conjured up from the mind of Moses Schönfinkel.

## Who Was Moses Schönfinkel?

So who was this person, who managed so long ago to see so far?

The complete known published output of Moses Schönfinkel consists of just two papers: his 1924 "On the Building Blocks of Mathematical Logic", and another, 31-page paper from 1927, coauthored with Paul Bernays, entitled "Zum *Entscheidungsproblem* der mathematischen Logik" ("On the Decision Problem of Mathematical Logic").

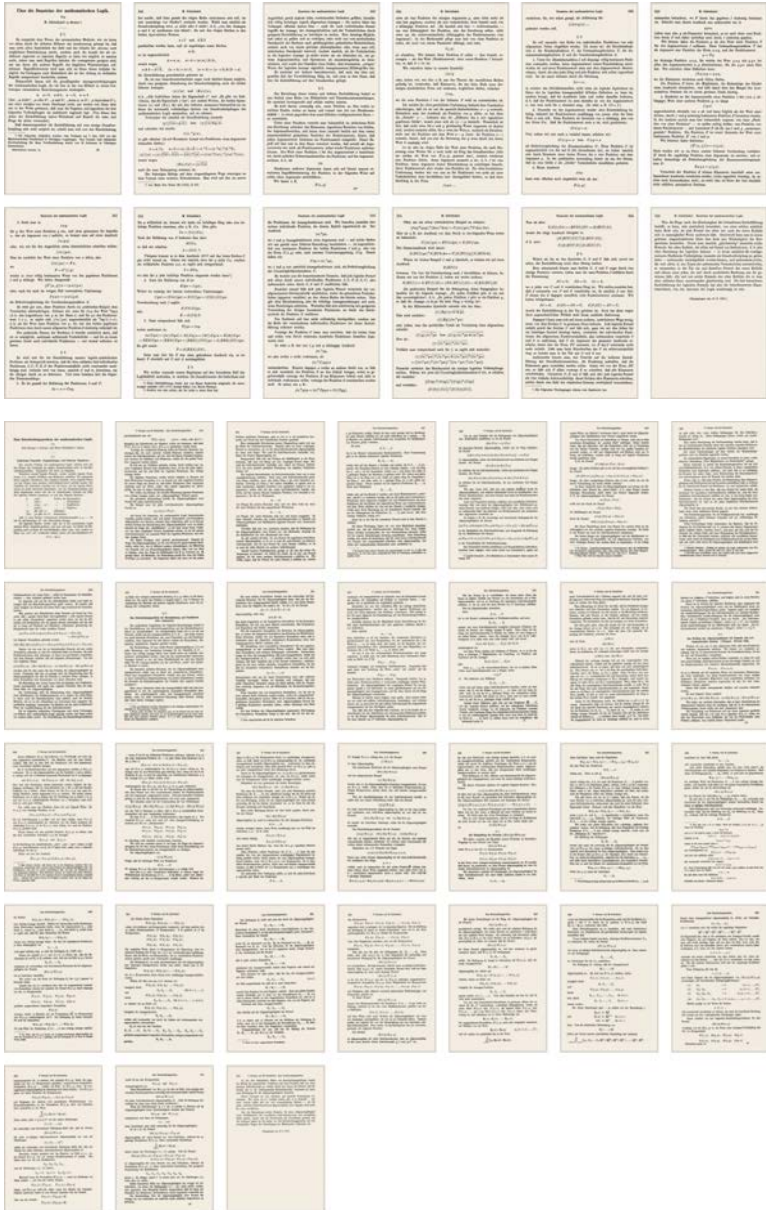

And somehow Schönfinkel has always been in the shadows—appearing at best only as a kind of footnote to a footnote. Turing machines have taken the limelight as models of computation—with combinators, hard to understand as they are, being mentioned at most only in obscure footnotes. And even within the study of combinators—often called "combinatory logic"—even as *S* and *K* have remained ubiquitous, Schönfinkel's invention of them typically garners at most a footnote.

About Schönfinkel as a person, three things are commonly said. First, that he was somehow connected with the mathematician David Hilbert in Göttingen. Second, that he spent time in a psychiatric institution. And third, that he died in poverty in Moscow, probably around 1940 or 1942.

But of course there has to be more to the story. And in recognition of the centenary of Schönfinkel's announcement of combinators, I decided to try to see what I could find out.

I don't think I've got all the answers. But it's been an interesting, if at times unsettling, trek through the Europe—and mathematics—of a century or so ago. And at the end of it I feel I've come to know and understand at least a little more about the triumph and tragedy of Moses Schönfinkel.

## The Beginning of the Story

It's a strange and sad resonance with Moses Schönfinkel's life... but there's a 1953 song by Tom Lehrer about plagiarism in mathematics—where the protagonist explains his chain of intellectual theft: "I have a friend in Minsk/Who has a friend in Pinsk/Whose friend in Omsk"... "/Whose friend somehow/Is solving now/The problem in Dnepropetrovsk". Well, Dnepropetrovsk is where Moses Schönfinkel was born.

Except, confusingly, at the time it was called (after Catherine the Great or maybe her namesake saint) Ekaterinoslav (Екатериносла́в)—and it's now called Dnipro. It's one of the larger cities in Ukraine, roughly in the center of the country, about 250 miles down the river Dnieper from Kiev. And at the time when Schönfinkel was born, Ukraine was part of the Russian Empire.

So what traces are there of Moses Schönfinkel in Ekaterinoslav (AKA Dnipro) today? 132 years later it wasn't so easy to find (especially during a pandemic)... but here's a record of his birth: a certificate from the Ekaterinoslav Public Rabbi stating that entry 272 of the Birth Register for Jews from 1888 records that on September 7, 1888, a son Moses was born to the Ekaterinoslav citizen Ilya Schönfinkel and his wife Masha:

This seems straightforward enough. But immediately there's a subtlety. When exactly was Moses Schönfinkel born? What is that date? At the time the Russian Empire—which had the Russian Orthodox Church, which eschewed Pope Gregory's 1582 revision of the calendar—was still using the Julian calendar introduced by Julius Caesar. (The calendar was switched in 1918 after the Russian Revolution, although the Orthodox Church plans to go on celebrating Christmas on January 7 until 2100.) So to know a correct modern (i.e. Gregorian calendar) date of birth we have to do a conversion. And from this we'd conclude that Moses Schönfinkel was born on September 19, 1888.

But it turns out that's not the end of the story. There are several other documents associated with Schönfinkel's college years that also list his date of birth as September 7, 1888. But the state archives of the Dnepropetrovsk region contain the actual, original register from the synagogue in Ekaterinoslav. And here's entry 272—and it records the birth of Moses Schönfinkel, but on September 17, not September 7:

So the official certificate is wrong! Someone left a digit out. And there's a check: the Birth Register also gives the date in the Jewish calendar: 24 Tishrei–which for 1888 is the Julian date September 17. So converting to modern Gregorian form, the correct date of birth for Moses Schönfinkel is September 29, 1888.

OK, now what about his name? In Russian it's given as Моисей Шейнфинкель (or, including the patronymic, with the most common transliteration from Hebrew, Моисей Эльевич Шейнфинкель). But how should his last name be transliterated? Well, there are several possibilities. We're using Schönfinkel—but other possibilities are Sheinfinkel and Sheynfinkel—and these show up almost randomly in different documents.

What else can we learn from Moses Schönfinkel's "birth certificate"? Well, it describes his father Эльева (Ilya) as an Ekaterinoslav мещанина. But what is that word? It's often translated "bourgeoisie", but seems to have basically meant "middle-class city dweller". And in other documents from the time, Ilya Schönfinkel is described as a "merchant of the 2nd guild" (i.e. not the "top 5%" 1st guild, nor the lower 3rd guild).

Apparently, however, his fortunes improved. The 1905 "Index of Active Enterprises Incorporated in the [Russian] Empire" lists him as a "merchant of the 1st guild" and records that in 1894 he co-founded the company of "Lurie & Sheinfinkel" (with a paid-in capital of 10,000 rubles, or about $150k today) that was engaged in the grocery trade:

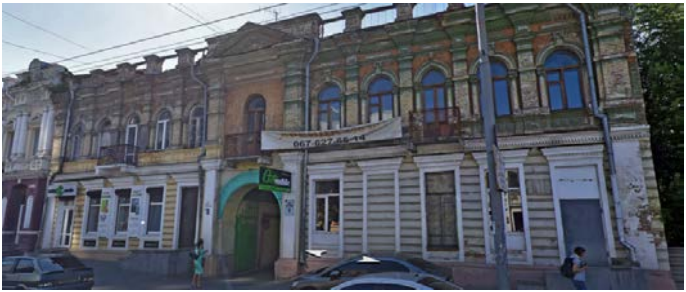

Lurie & Sheinfinkel seems to have had multiple wine and grocery stores. Between 1901 and 1904 its "store #2" was next to a homeopathic pharmacy in a building that probably looked at the time much like it does today:

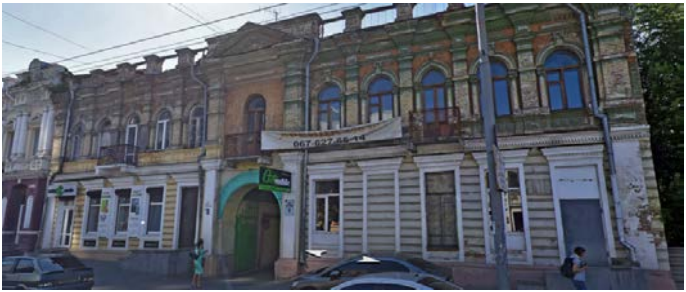

And for store #1 there are actually contemporary photographs (note the -инкель for the end of "Schönfinkel" visible on the bottom left; this particular building was destroyed in World War II):

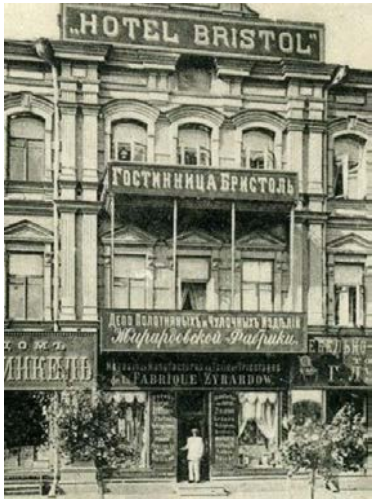

There seems to have been a close connection between the Schönfinkels and the Luries—who were a prominent Ekaterinoslav family involved in a variety of enterprises. Moses Schönfinkel's mother Maria (Masha) was originally a Lurie (actually, she was one of the 8 siblings of Ilya Schönfinkel's business partner Aron Lurie). Ilya Schönfinkel is listed from 1894 to 1897 as "treasurer of the Lurie Synagogue". And in 1906 Moses Schönfinkel listed his mailing address in Ekaterinoslav as Lurie House, Ostrozhnaya Square. (By 1906 that square sported an upscale park—though a century earlier it had housed a prison that was referenced in a poem by Pushkin. Now it's the site of an opera house.)

Accounts of Schönfinkel sometimes describe him as coming from a "village in Ukraine". In actuality, at the turn of the twentieth century Ekaterinoslav was a bustling metropolis, that for example had just become the third city in the whole Russian Empire to have electric trams. Schönfinkel's family also seems to have been quite well to do. Some pictures of Ekaterinoslav from the time give a sense of the environment (this building was actually the site of a Lurie candy factory):

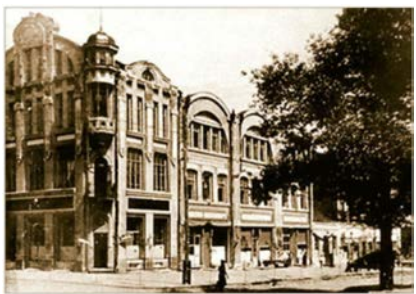

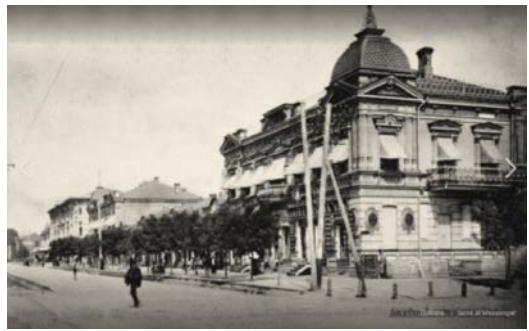

As the name "Moses" might suggest, Moses Schönfinkel was Jewish, and at the time he was born there was a large Jewish population in the southern part of Ukraine. Many Jews had come to Ekaterinoslav from Moscow, and in fact 40% of the whole population of the town was identified as Jewish.

Moses Schönfinkel went to the main high school in town (the "Ekaterinoslav classical gymnasium")—and graduated in 1906, shortly before turning 18. Here's his diploma:

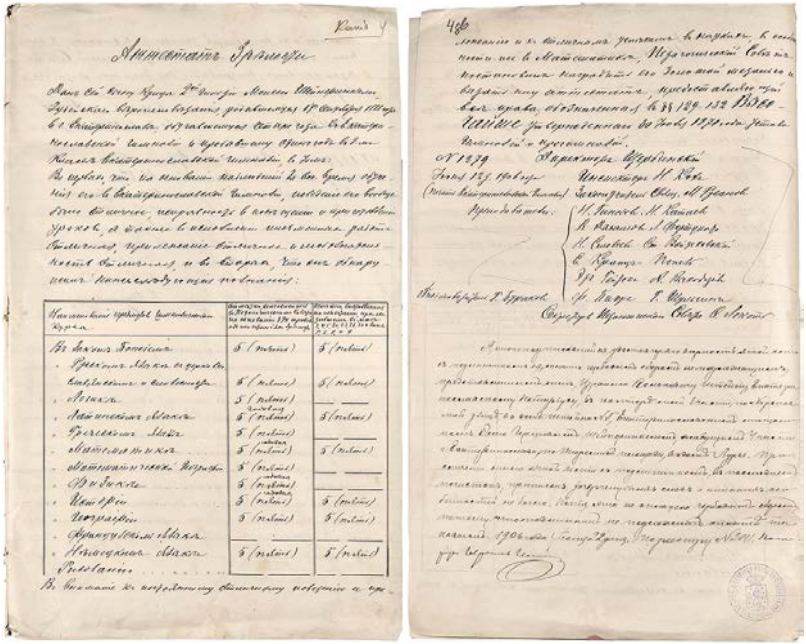

The diploma shows that he got 5/5 in all subjects—the subjects being theology, Russian, logic, Latin, Greek, mathematics, geodesy ("mathematical geography"), physics, history, geography, French, German and drawing. So, yes, he did well in high school. And in fact the diploma goes on to say: "In view of his excellent behavior and diligence and excellent success in the sciences, especially in mathematics, the Pedagogical Council decided to award him the Gold Medal…"

## Going to College in Odessa

Having graduated from high school, Moses Schönfinkel wanted to go ("for purely family reasons", he said) to the University of Kiev. But being told that Ekaterinoslav was in the wrong district for that, he instead asked to enroll at Novorossiysk University in Odessa. He wrote a letter—in rather neat handwriting—to unscramble a bureaucratic issue, giving various excuses along the way:

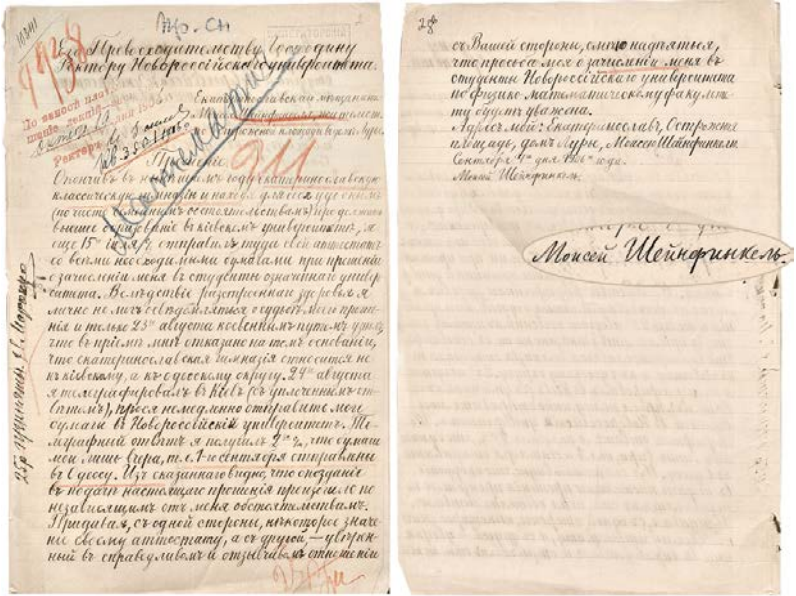

But in the fall of 1906, there he was: a student in the Faculty of Physics and Mathematics Faculty of Novorossiysk University, in the rather upscale and cosmopolitan town of Odessa, on the Black Sea.

The Imperial Novorossiya University, as it was then officially called, had been created out of an earlier institution by Tsar Alexander II in 1865. It was a distinguished university, with for example Dmitri Mendeleev (of periodic table fame) having taught there. In Soviet times it would be renamed after the discoverer of macrophages, Élie Metchnikoff (who worked there). Nowadays it is usually known as Odessa University. And conveniently, it has maintained its archives well—so that, still there, 114 years later, is Moses Schönfinkel's student file:

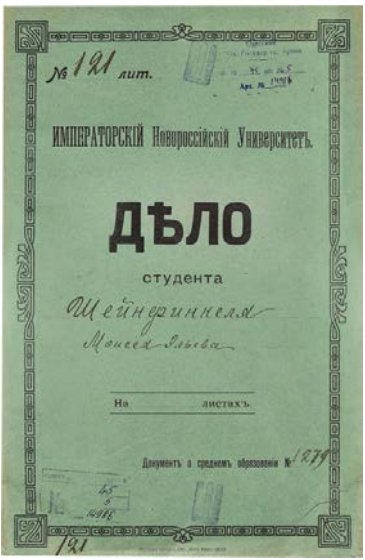

It's amazing how "modern" a lot of what's in it seems. First, there are documents Moses Schönfinkel sent so he could register (confirming them by telegram on September 1, 1906). There's his high-school diploma and birth certificate—and there's a document from the Ekaterinoslav City Council certifying his "citizen rank" (see above). The cover sheet also records a couple of other documents, one of which is presumably some kind of deferment of military service.

And then in the file there are two "photo cards" giving us pictures of the young Moses Schönfinkel, wearing the uniform of the Imperial Russian Army:

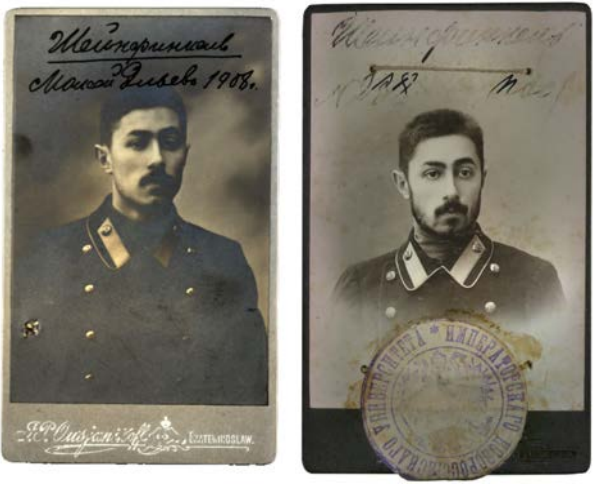

(These pictures actually seem to come from 1908; the style of uniform was a standard one issued after 1907; the [presumably] white collar tabs indicate the 3rd regiment of whatever division he was assigned to.)

Nowadays it would all be online, but in his physical file there is a "lecture book" listing courses (yes, every document is numbered, to correspond to a line in a central ledger):

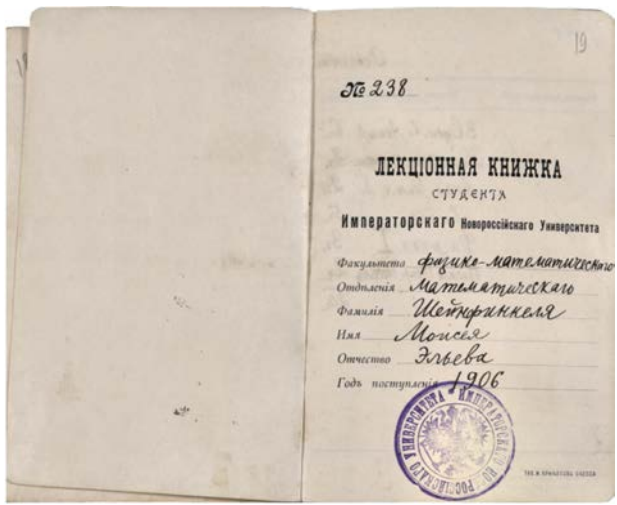

Here are the courses Moses Schönfinkel took in his first semester in college (fall 1906):

Introduction to Analysis (6 hrs), Introduction to Determinant Theory (2 hrs), Analytical Geometry 1 (2 hrs), Chemistry (5 hrs), Physics 1 (3 hrs), Elementary Number Theory (2 hrs): a total of 20 hours. Here's the bill for these courses: pretty good value at 1 ruble per course-hour, or a total of 20 rubles, which is about $300 today:

Subsequent semesters list many very familiar courses: Differential Calculus, Integrals (parts 1 and 2), and Higher Algebra, as well as "Calculus of Probabilities" (presumably probability theory) and "Determinant Theory" (essentially differently branded "linear algebra"). There are some "distribution" courses, like Astronomy (and Spherical Astronomy) and Physical Geography (or is that Geodesy?). And by 1908, there are also courses like Functions of a Complex Variable, Integro-Differential Equations (yeah, differential equations definitely pulled ahead of integral equations over the past century), Calculus of Variations and Infinite Series. And—perhaps presaging Schönfinkel's next life move—another course that makes an appearance in 1908 is German (and it's Schönfinkel's only non-science course during his whole university career).

In Schönfinkel's "lecture book" many of the courses also have names of professors listed. For example, there's "Kagan", who's listed as teaching Foundations of Geometry (as well as Higher Algebra, Determinant Theory and Integro-Differential Equations). That's Benjamin Kagan, who was then a young lecturer, but would later become a leader in differential geometry in Moscow—and also someone who studied the axiomatic foundations of geometry (as well as writing about the somewhat tragic life of Lobachevsky).

Another professor—listed as teaching Schönfinkel Introduction to Analysis and Theory of Algebraic Equation Solving—is "Shatunovsky". And (at least according to Shatunovsky's later student Sofya Yanoskaya, of whom we'll hear more later), Samuil Shatunovsky was basically Schönfinkel's undergraduate advisor.

Shatunovsky had been the 9th child of a poor Jewish family (actually) from a village in Ukraine. He was never able to enroll at a university, but for some years did manage to go to lectures by people around Pafnuty Chebyshev in Saint Petersburg. For quite a few years he then made a living as an itinerant math tutor (notably in Ekaterinoslav) but papers he wrote were eventually noticed by people at the university in Odessa, and, finally, in 1905, at the age of 46, he ended up as a lecturer at the university—where the following year he taught Schönfinkel.

Shatunovsky (who stayed in Odessa until his death in 1929) was apparently an energetic but precise lecturer. He seems to have been quite axiomatically oriented, creating axiomatic systems for geometry, algebraic fields, and notably, for order relations. (He was also quite a constructivist, opposed to the indiscriminate use of the Law of Excluded Middle.) The lectures from his Introduction to Analysis course (which Schönfinkel took in 1906) were published in 1923 (by the local publishing company Mathesis in which he and Kagan were involved).

Another of Schönfinkel's professors (from whom he took Differential Calculus and "Calculus of Probabilities") was a certain Ivan (or Jan) Śleszyński, who had worked with Karl Weierstrass on things like continued fractions, but by 1906 was in his early 50s and increasingly transitioning to working on logic. In 1911 he moved to Poland, where he sowed some of the seeds for the Polish school of mathematical logic, in 1923 writing a book called *On the Significance of Logic for Mathematics* (notably with no mention of Schönfinkel), and in 1925 one on proof theory.

It's not clear how much mathematical logic Moses Schönfinkel picked up in college, but in any case, in 1910, he was ready to graduate. Here's his final student ID (what are those pieces of string for?):

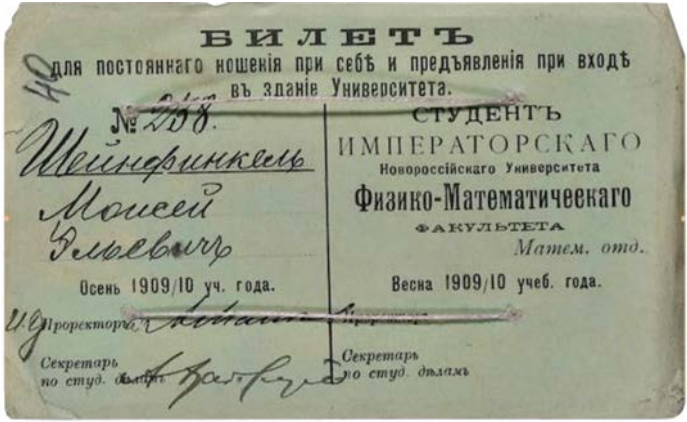

There's a certificate confirming that on April 6, 1910, Moses Schönfinkel had no books that needed returning to the library. And he sent a letter asking to graduate (with slightly-less-neat handwriting than in 1906):

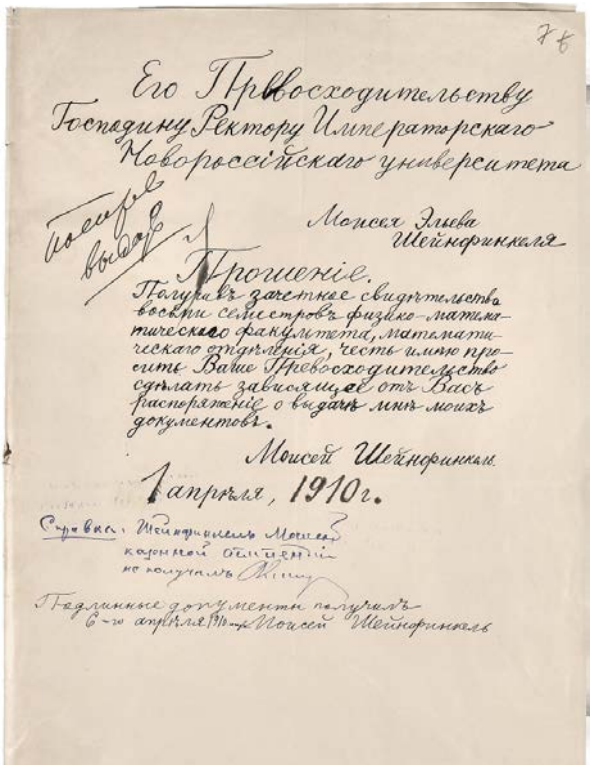

The letter closes with his signature (Моисей Шейнфинкель):

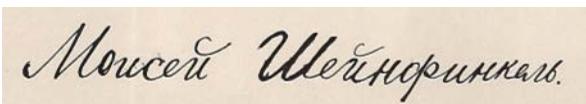

# Göttingen, Center of the Mathematical Universe

After Moses Schönfinkel graduated college in 1910 he probably went into four years of military service (perhaps as an engineer) in the Russian Imperial Army. World War I began on July 28, 1914—and Russia mobilized on July 30. But in one of his few pieces of good luck Moses Schönfinkel was not called up, having arrived in Göttingen, Germany, on June 1, 1914 (just four weeks before the event that would trigger World War I), to study mathematics.

Göttingen was at the time a top place for mathematics. In fact, it was sufficiently much of a "math town" that around that time postcards of local mathematicians were for sale there. And the biggest star was David Hilbert—which is who Schönfinkel went to Göttingen hoping to work with.

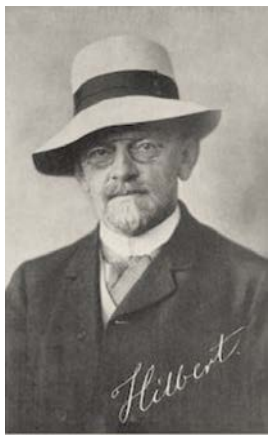

Hilbert had grown up in Prussia and started his career in Königsberg. His big break came in 1888 at age 26 when he got a major result in representation theory (then called "invariant theory")—using then-shocking non-constructive techniques. And it was soon after this that Felix Klein recruited Hilbert to Göttingen—where he remained for the rest of his life.

In 1900 Hilbert gave his famous address to the International Congress of Mathematicians where he first listed his (ultimately 23) problems that he thought should be important in the future of mathematics. Almost all the problems are what anyone would call "mathematical". But problem 6 has always stuck out for me: "Mathematical Treatment of the Axioms of Physics": Hilbert somehow wanted to axiomatize physics as Euclid had axiomatized geometry. And he didn't just talk about this; he spent nearly 20 years working on it. He brought in physicists to teach him, and he worked on things like gravitation theory ("Einstein–Hilbert action") and kinetic theory—and wanted for example to derive the existence of the electron from something like Maxwell's equations. (He was particularly interested in the way atomistic processes limit to continua—a problem that I now believe is deeply connected to computational irreducibility, in effect implying another appearance of undecidability, like in Hilbert's 1st, 2nd and 10th problems.)

Hilbert seemed to feel that physics was a crucial source of raw material for mathematics. But yet he developed a whole program of research based on doing mathematics in a completely formalistic way—where one just writes down axioms and somehow "mechanically" generates all true theorems from them. (He seems to have drawn some distinction between "merely mathematical" questions, and questions about physics, apparently noting—in a certain resonance with my life's work—that in the latter case "the physicist has the great calculating machine, Nature".)

In 1899 Hilbert had written down more precise and formal axioms for Euclid's geometry, and he wanted to go on and figure out how to formulate other areas of mathematics in this kind of axiomatic way. But for more than a decade he seems to have spent most of his time on physics—finally returning to questions about the foundations of mathematics around 1917, giving lectures about "logical calculus" in the winter session of 1920.

By 1920, World War I had come and gone, with comparatively little effect on mathematical life in Göttingen (the nearest battle was in Belgium 200 miles to the west). Hilbert was 58 years old, and had apparently lost quite a bit of his earlier energy (not least as a result of having contracted pernicious anemia [autoimmune vitamin B12 deficiency], whose cure was found only a few years later). But Hilbert was still a celebrity around Göttingen, and generating mathematical excitement. (Among "celebrity gossip" mentioned in a letter home by young Russian topologist Pavel Urysohn is that Hilbert was a huge fan of the gramophone, and that even at his advanced age, in the summer, he would sit in a tree to study.)

I have been able to find out almost nothing about Schönfinkel's interaction with Hilbert. However, from April to August 1920 Hilbert gave weekly lectures entitled "Problems of Mathematical Logic" which summarized the standard formalism of the field—and the official notes for those lectures were put together by Moses Schönfinkel and Paul Bernays (the "N" initial for Schönfinkel is a typo):

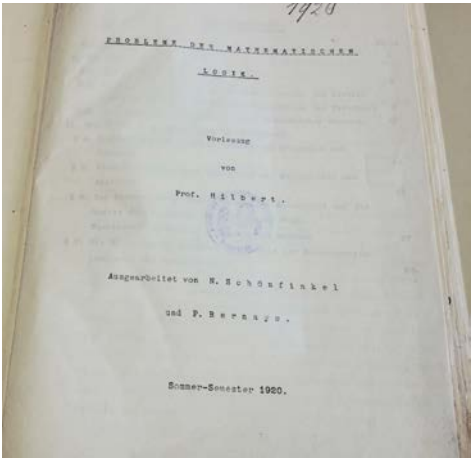

*Photograph by Cem Bozşahin*

A few months after these lectures came, at least from our perspective today, the highlight of Schönfinkel's time in Göttingen: the talk he gave on December 7, 1920. The venue was the weekly meeting of the Göttingen Mathematics Society, held at 6pm on Tuesdays. The society wasn't officially part of the university, but it met in the same university "Auditorium Building" that at the time housed the math institute:

The talks at the Göttingen Mathematics Society were listed in the *Annual Report of the German Mathematicians Association*:

There's quite a lineup. November 9, Ludwig Neder (student of Edmund Landau): "Trigonometric Series". November 16, Erich Bessel-Hagen (student of Carathéodory): "Discontinuous Solutions of Variational Problems". November 23, Carl Runge (of Runge–Kutta fame, then a Göttingen professor): "American Work on Star Clusters in the Milky Way". November 30 Gottfried Rückle (assistant of van der Waals): "Explanations of Natural Laws Using a Statistical Mechanics Basis". And then: December 7: Moses Schönfinkel, "Elements of Logic".

The next week, December 14, Paul Bernays, who worked with Hilbert and interacted with Schönfinkel, spoke about "Probability, the Arrow of Time and Causality" (yes, there was still a lot of interest around Hilbert in the foundations of physics). January 10+11, Joseph Petzoldt (philosopher of science): "The Epistemological Basis of Special and General Relativity". January 25, Emmy Noether (of Noether's theorem fame): "Elementary Divisors and General Ideal Theory". February 1+8, Richard Courant (of PDE etc. fame) & Paul Bernays: "About the New Arithmetic Theories of Weyl and Brouwer". February 22, David Hilbert: "On a New Basis for the Meaning of a Number" (yes, that's foundations of math).

What in detail happened at Schönfinkel's talk, or as a result of it? We don't know. But he seems to have been close enough to Hilbert that just over a year later he was in a picture taken for David Hilbert's 60th birthday on January 23, 1922:

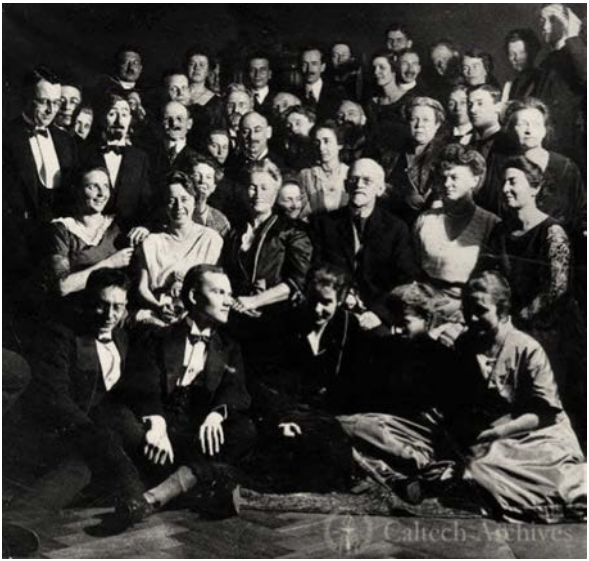

There are all sorts of well-known mathematicians in the picture (Richard Courant, Hermann Minkowski, Edmund Landau, ...) as well as some physicists (Peter Debye, Theodore von Kármán, Ludwig Prandtl, ...). And there near the top left is Moses Schönfinkel, sporting a somewhat surprised expression.

For his 60th birthday Hilbert was given a photo album—with 44 pages of pictures of altogether about 200 mathematicians (and physicists). And there on page 22 is Moses Schönfinkel:

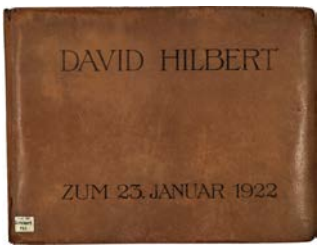

*Göttingen University, Cod. Ms. D. Hilbert 754*

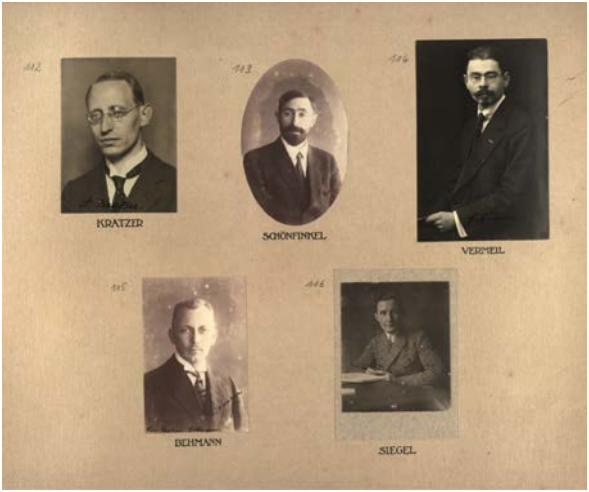

*Göttingen University, Cod. Ms. D. Hilbert 754, Bl. 22*

Who are the other people on the page with him?  Adolf Kratzer (1893–1983) was a student of Arnold Sommerfeld, and a "physics assistant" to Hilbert.  Hermann Vermeil (1889–1959) was an assistant to Hermann Weyl, who worked on differential geometry for general relativity. Heinrich Behmann (1891–1970) was a student of Hilbert and worked on mathematical logic, and we'll encounter him again later.  Finally,  Carl Ludwig Siegel (1896–1981) had been a student of Landau and would become a well-known number theorist.

# Problems Are Brewing

There's a lot that's still mysterious about Moses Schönfinkel's time in Göttingen. But we have one (undated) letter written by Nathan Schönfinkel, Moses's younger brother,  presumably in 1921 or 1922 (yes, he romanizes his name "Scheinfinkel" rather than "Schönfinkel"):

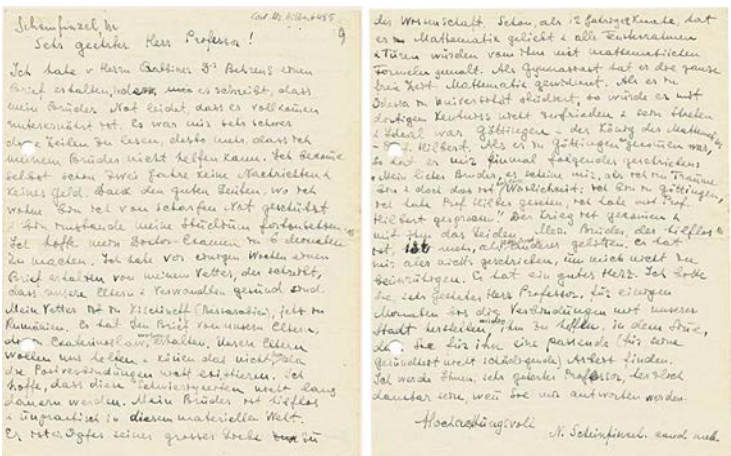

*Göttingen University, Cod. Ms. D. Hilbert 455: 9*

Dear Professor!

I received a letter from Rabbi Dr. Behrens in which he wrote that my brother was in need, that he was completely malnourished. It was very difficult for me to read these lines, even more so because I cannot help my brother. I haven't received any messages or money myself for two years. Thanks to the good people where I live, I am protected from severe hardship. I am able to continue my studies. I hope to finish my PhD in 6 months. A few weeks ago I received a letter from my cousin stating that our parents and relatives are healthy. My cousin is in Kishinev (Bessarabia), now in Romania. He received the letter from our parents who live in Ekaterinoslav. Our parents want to help us but cannot do so because the postal connections are nonexistent. I hope these difficulties will not last long. My brother is helpless and impractical in this material world. He is a victim of his great love for science. Even as a 12 year old boy he loved mathematics, and all window frames and doors were painted with mathematical formulas by him. As a high school student, he devoted all his free time to mathematics. When he was studying at the university in Odessa, he was not satisfied with the knowledge there, and his striving and ideal was Göttingen and the king of mathematics, Prof. Hilbert. When he was accepted in Göttingen, he once wrote to me the following: "My dear brother, it seems to me as if I am dreaming but this is reality: I am in Göttingen, I saw Prof. Hilbert, I spoke to Prof. Hilbert." The war came and with it suffering. My brother, who is helpless, has suffered more than anyone else. But he did not write to me so as not to worry me. He has a good heart. I ask you, dear Professor, for a few months until the connections with our city are established, to help him by finding a suitable (not harmful to his health) job for him. I will be very grateful to you, dear Professor, if you will answer me.

Sincerely.

N. Scheinfinkel

We'll talk more about Nathan Schönfinkel later. But suffice it to say here that when he wrote the letter he was a physiology graduate student at the University of Bern—and he would get his PhD in 1922, and later became a professor. But the letter he wrote is probably our best single surviving source of information about the situation and personality of Moses Schönfinkel. Obviously he was a serious math enthusiast from a young age. And the letter implies that he'd wanted to work with Hilbert for some time (presumably hence the German classes in college).

It also implies that he was financially supported in Göttingen by his parents—until this was disrupted by World War I. (And we learn that his parents were OK in the Russian Revolution.) (By the way, the rabbi mentioned is probably a certain Siegfried Behrens, who left Göttingen in 1922.)

There's no record of any reply to Nathan Schönfinkel's letter from Hilbert. But at least by the time of Hilbert's 60th birthday in 1922 Moses Schönfinkel was (as we saw above) enough in the inner circle to be invited to the birthday party.

What else is there in the university archives in Göttingen about Moses Schönfinkel? There's just one document, but it's very telling:

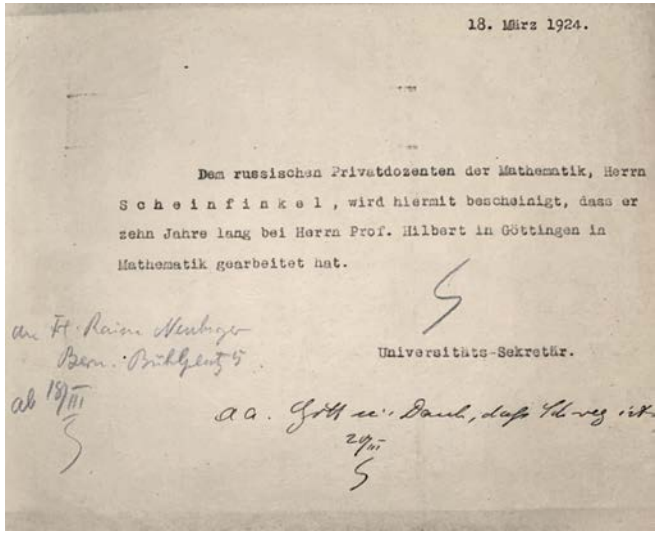

It's dated 18 March 1924. And it's a carbon copy of a reference for Schönfinkel. It's rather cold and formal, and reads:

> "The Russian privatdozent [private lecturer] in mathematics, Mr. Scheinfinkel, is hereby certified to have worked in mathematics for ten years with Prof. Hilbert in Göttingen."

It's signed (with a stylized "S") by the "University Secretary", a certain Ludwig Gossmann, who we'll be talking about later. And it's being sent to Ms. Raissa Neuburger, at Bühlplatz 5, Bern. That address is where the Physiology Institute at the University of Bern is now, and also was in 1924. And Raissa Neuberger either was then, or soon would become, Nathan Schönfinkel's wife.

But there's one more thing, handwritten in black ink at the bottom of the document. Dated March 20, it's another note from the University Secretary. It's annotated "a.a.", i.e. *ad acta*— for the records. And in German it reads:

> *Gott sei Dank, dass Sch weg ist*

which translates in English as:

> Thank goodness Sch is gone

Hmm. So for some reason at least the university secretary was happy to see Schönfinkel go. (Or perhaps it was a German 1920s version of an HR notation: "not eligible for rehire".) But let's analyze this document in a little more detail. It says Schönfinkel worked with Hilbert for 10 years. That agrees with him having arrived in Göttingen in 1914 (which is a date we know for other reasons, as we'll see below).

But now there's a mystery. The reference describes Schönfinkel as a "*privatdozent*". That's a definite position at a German university, with definite rules, that in 1924 one would expect to have been rigidly enforced. The basic career track was (and largely still is): first, spend 2–5 years getting a PhD. Then perhaps get recruited for a professorship, or if not, continue doing research, and write a habilitation, after which the university may issue what amounts to an official government "license to teach", making someone a privatdozent, able to give lectures. Being a privatdozent wasn't as such a paid gig. But it could be combined with a job like being an assistant to a professor—or something outside the university, like tutoring, teaching high school or working at a company.

So if Schönfinkel was a privatdozent in 1924, where is the record of his PhD, or his habilitation? To get a PhD required "formally publishing" a thesis, and printing (as in, on a printing press) at least 20 or so copies of the thesis. A habilitation was typically a substantial, published research paper. But there's absolutely no record of any of these things for Schönfinkel. And that's very surprising. Because there are detailed records for other people (like Paul Bernays) who were around at the time, and were indeed privatdozents.

And what's more the *Annual Report of the German Mathematicians Association*—which listed Schönfinkel's 1920 talk—seems to have listed mathematical goings-on in meticulous detail. Who gave what talk. Who wrote what paper. And most definitely who got a PhD, did a habilitation or became a privatdozent. (And becoming a privatdozent also required an action of the university senate, which was carefully recorded.) But going through all the annual reports of the German Mathematicians Association we find only four mentions of Schönfinkel. There's his 1920 talk, and also a 1921 talk with Paul Bernays that we'll discuss later. There's the publication of his papers in 1924 and 1927. And there's a single other entry, which says that on November 4, 1924, Richard Courant gave a report to the Göttingen Mathematical Society about a conference in Innsbruck, where Heinrich Behmann reported on "published work by M. Schönfinkel". (It describes the work as follows: "It is a continuation of Sheffer's [1913] idea of replacing the elementary operations of symbolic logic with a single one. By means of a certain function calculus, all logical statements (including the mathematical ones) are represented by three basic signs alone.")

So, it seems, the university secretary wasn't telling it straight. Schönfinkel might have worked with Hilbert for 10 years. But he wasn't a privatdozent. And actually it doesn't seem as if he had any "official status" at all.

So how do we even know that Schönfinkel was in Göttingen from 1914 to 1924? Well, he was Russian, and so in Germany he was an "alien", and as such he was required to register his address with the local police (no doubt even more so from 1914 to 1918 when Germany was, after all, at war with Russia). And the remarkable thing is that even after all these years, Schönfinkel's registration card is still right there in the municipal archives of the city of Göttingen:

*Stadtarchiv Göttingen, Meldekartei*

So that means we have all Schönfinkel's addresses during his time in Göttingen. Of course, there are confusions. There's yet another birthdate for Schönfinkel: September 4, 1889. Wrong year. Perhaps a wrongly done correction from the Julian calendar. Perhaps "adjusted" for some reason of military service obligations. But, in any case, the document says that Moses Schönfinkel from Ekaterinoslav arrived in Göttingen on June 1, 1914, and started living at 6 Lindenstraße (now Felix-Klein-Strasse).

He moved pretty often (11 times in 10 years), not at particularly systematic times of year. It's not clear exactly what the setup was in all these places, but at least at the end (and in another document) it lists addresses and "with Frau ….", presumably indicating that he was renting a room in someone's house.

Where were all those addresses? Well, here's a map of Göttingen circa 1920, with all of them plotted (along with a red "M" for the location of the math institute):

*Stadtarchiv Göttingen, D 2, V a 62*

The last item on the registration card says that on March 18, 1924, he departed Göttingen, and went to Moscow. And the note on the copy of the reference saying "thank goodness [he's] gone" is dated March 20, so that all ties together.

But let's come back to the reference. Who was this "University Secretary" who seems to have made up the claim that Schönfinkel was a privatdozent? It was fairly easy to find out that his name was Ludwig Gossmann. But the big surprise was to find out that the university archives in Göttingen have nearly 500 pages about him—primarily in connection with a "criminal investigation".

Here's the story. Ludwig Gossmann was born in 1878 (so he was 10 years older than Schönfinkel). He grew up in Göttingen, where his father was a janitor at the university. He finished high school but didn't go to college and started working for the local government. Then in 1906 (at age 28) he was hired by the university as its "secretary".

The position of "university secretary" was a high-level one. It reported directly to the vice-rector of the university, and was responsible for "general administrative matters" for the university, including, notably, the supervision of international students (of whom there were many, Schönfinkel being one). Ludwig Gossmann held the position of university secretary for 27 years—even while the university had a different rector (normally a distinguished academic) every year.

But Mr. Gossmann also had a sideline: he was involved in real estate. In the 1910s he started building houses (borrowing money from, among others, various university professors). And by the 1920s he had significant real estate holdings—and a business where he rented to international visitors and students at the university.

Years went by. But then, on January 24, 1933, the newspaper headline announced: "Sensational arrest: senior university official Gossmann arrested on suspicion of treason—communist revolution material [*Zersetzungsschrift*] confiscated from his apartment". It was said that perhaps it was a setup, and that he'd been targeted because he was gay (though, a year earlier, at age 54, he did marry a woman named Elfriede).

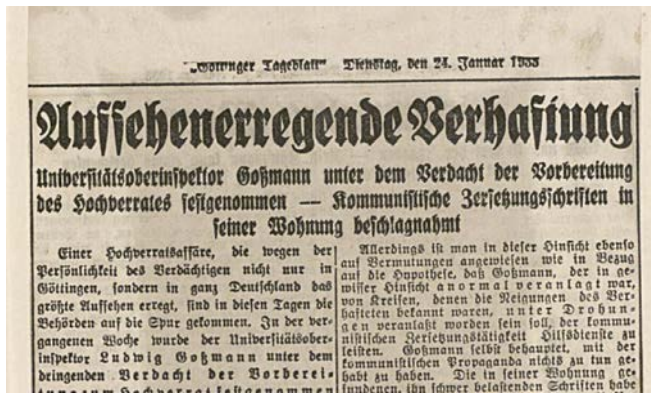

*Göttingen University, Kur 3730, Sek 356 2*

This was a bad time to be accused of being a communist (Hitler would become chancellor less than a week later, on January 30, 1933, in part propelled by fears of communism). Gossmann was taken to Hanover "for questioning", but was then allowed back to Göttingen "under house arrest". He'd had health problems for several years, and died of a heart attack on February 24, 1933.

But none of this really helps us understand why Gossmann would go out on a limb to falsify the reference for Schönfinkel. We can't specifically find an address match, but perhaps Schönfinkel had at least at some point been a tenant of Gossmann's. Perhaps he still owed rent. Perhaps he was just difficult in dealing with the university administration. It's not clear. It's also not clear why the reference Gossmann wrote was sent to Schönfinkel's brother in Bern, even though Schönfinkel himself was going to Moscow. Or why it wasn't just handed to Schönfinkel before he left Göttingen.

## The 1924 Paper

Whatever was going on with Schönfinkel in Göttingen in 1924, we know one thing for sure: it was then that he published his remarkable paper about what are now called combinators. Let's talk in a bit more detail about the paper—though the technicalities I'm discussing elsewhere.

First, there's some timing. At the end of the paper, it says it was received by the journal on March 15, 1924, i.e. just three days before the date of Ludwig Gossmann's reference for Schönfinkel. And then at the top of the paper, there's something else: under Schönfinkel's name it says "in Moskau", i.e. at least as far as the journal was concerned, Schönfinkel was in Moscow, Russia, at the time the article was published:

Über die Bausteine der mathematischen Logik.

Von

M. Schönfinkel in Moskau[1]).

There's also a footnote on the first page of the paper:

[1]) Die folgenden Gedanken wurden vom Verfasser am 7. Dez. 1920 vor der Mathematischen Gesellschaft in Göttingen vorgetragen. Ihre formale und stilistische Durcharbeitung für diese Veröffentlichung wurde von H. Behmann in Göttingen übernommen.

"The following thoughts were presented by the author to the Mathematical Society in Göttingen on December 7, 1920. Their formal and stylistic processing for this publication was done by H. Behmann in Göttingen."

The paper itself is written in a nice, clear and mathematically mature way. Its big result (as I've discussed elsewhere) is the introduction of what would later be called combinators: two abstract constructs from which arbitrary functions and computations can be built up. Schönfinkel names one of them *S*, after the German word "*Schmelzen*" for "fusion". The other has become known as *K*, although Schönfinkel calls it *C*, even though the German word for "constancy" (which is what would naturally describe it) is "*Konstantheit*", which starts with a K.

The paper ends with three paragraphs, footnoted with "The considerations that follow are the editor's" (i.e. Behmann's). They're not as clear as the rest of the paper, and contain a confused mistake.

The main part of the paper is "just math" (or computation, or whatever). But here's the page where *S* and *K* (called *C* here) are first used:

And now there's something more people-oriented: a footnote to the combinator equation *I = SCC* saying "This reduction was communicated to me by Mr. Boskowitz; some time before that, Mr. Bernays had called the somewhat less simple one (*SC*)(*CC*) to my attention." In other words, even if nothing else, Schönfinkel had talked to Boskowitz and Bernays about what he was doing.

OK, so we've got three people—in addition to David Hilbert—somehow connected to Moses Schönfinkel.

Let's start with Heinrich Behmann—the person footnoted as "processing" Schönfinkel's paper for publication:

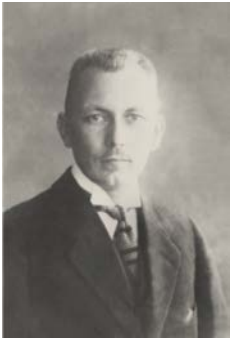

He was born in Bremen, Germany, in 1891, making him a couple of years younger than Schönfinkel. He arrived in Göttingen as a student in 1911, and by 1914 was giving a talk about Whitehead and Russell's *Principia Mathematica* (which had been published in 1910). When World War I started he volunteered for military service, and in 1915 he was wounded in action in Poland (receiving an Iron Cross)—but in 1916 he was back in Göttingen studying under Hilbert, and in 1918 he wrote his PhD thesis on "The Antinomy of the Transfinite Number and Its Resolution by the Theory of Russell and Whitehead" (i.e. using the idea of types to deal with paradoxes associated with infinity).

Behmann continued in the standard academic track (i.e. what Schönfinkel apparently didn't do)—and in 1921 he got his habilitation with the thesis "Contributions to the Algebra of Logic, in Particular to the *Entscheidungsproblem* [Decision Problem]". There'd been other decision problems discussed before, but Behmann said what he meant was a "procedure [giving] complete instructions for determining whether a [logical or mathematical] assertion is true or false by a deterministic calculation after finitely many steps". And, yes, Alan Turing's 1936 paper "On Computable Numbers, with an Application to the *Entscheidungsproblem*" was what finally established that the halting problem, and therefore the *Entscheidungsproblem*, was undecidable. Curiously, in principle, there should have been enough in Schönfinkel's paper that this could have been figured out back in 1921 if Behmann or others had been thinking about it in the right way (which might have been difficult before Gödel's work).

So what happened to Behmann? He continued to work on mathematical logic and the philosophy of mathematics. After his habilitation in 1921 he became a privatdozent at Göttingen (with a job as an assistant in the applied math institute), and then in 1925 got a professorship in Halle in applied math—though having been an active member of the Nazi Party since 1937, lost this professorship in 1945 and became a librarian. He died in 1970.

(By the way, even though in 1920 "PM" [*Principia Mathematica*] was hot—and Behmann was promoting it—Schönfinkel had what in my opinion was the good taste to not explicitly mention it in his paper, referring only to Hilbert's much-less-muddy ideas about the formalization of mathematics.)

OK, so what about Boskovitz, credited in the footnote with having discovered the classic combinator result *I = SKK*? That was Alfred Boskovitz, in 1920 a 23-year-old Jewish student at Göttingen, who came from Budapest, Hungary, and worked with Paul Bernays on set theory. Boskovitz is notable for having contributed far more corrections (nearly 200) to *Principia Mathematica* than anyone else, and being acknowledged (along with Behmann) in a footnote in the (1925–27) second edition. (This edition also gives a reference to Schönwinkel's [sic] paper at the end of a list of 14 "other contributions to mathematical logic" since the first edition.) In the mid-1920s Boskovitz returned to Budapest. In 1936 he wrote to Behmann that anti-Jewish sentiment there made him concerned for his safety. There's one more known communication from him in 1942, then no further trace.

The third person mentioned in Schönfinkel's paper is Paul Bernays, who ended up living a long and productive life, mostly in Switzerland. But we'll come to him later.

So where was Schönfinkel's paper published? It was in a journal called *Mathematische Annalen* (*Annals of Mathematics*)—probably the top math journal of the time. Here's its rather swank masthead, with quite a collection of famous names (including physicists like Einstein, Born and Sommerfeld):

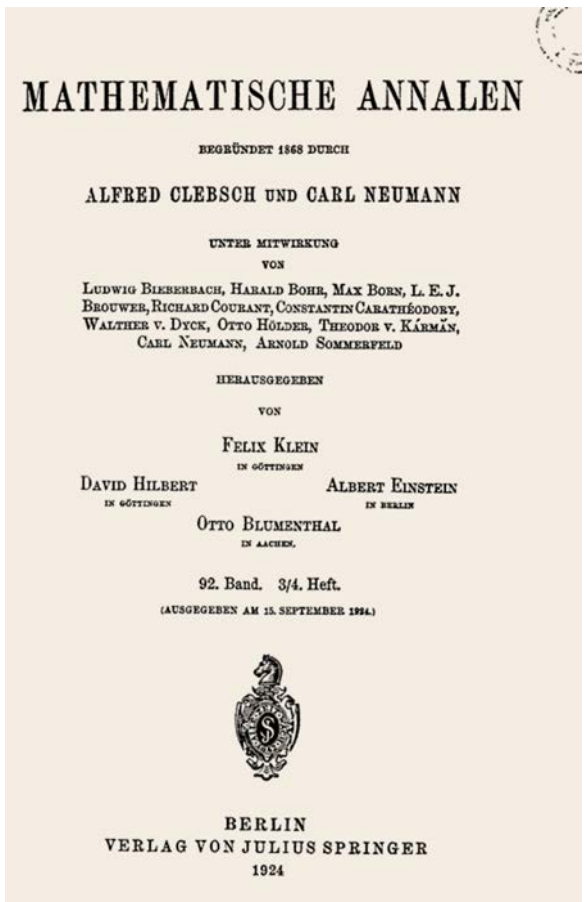

The "instructions to contributors" on the inside cover of each issue had a statement from the "Editorial Office" about not changing things at the proof stage because "according to a calculation they [cost] 6% of the price of a volume". The instructions then go on to tell people to submit papers to the editors—at their various home addresses (it seems David Hilbert lived just down the street from Felix Klein...):

Here's the complete table of contents for the volume in which Schönfinkel's paper appears:

There are a variety of famous names here. But particularly notable for our purposes are Aleksandr Khintchine (of Khinchin constant fame) and the topologists Pavel Alexandroff and Pavel Urysohn, who were all from Moscow State University, and who are all indicated, like Schönfinkel, as being "in Moscow".

There's a little bit of timing information here. Schönfinkel's paper was indicated as having been received by the journal on March 15, 1924. The "thank goodness [he's] gone [from Göttingen]" comment is dated March 20. Meanwhile, the actual issue of the journal with Schönfinkel's article (number 3 of 4) was published September 15, with table of contents:

But note the ominous † next to Urysohn's name. Turns out his fatal swimming accident was August 17, so—notwithstanding their admonitions—the journal must have added the † quite quickly at the proof stage.

## The "1927" Paper

Beyond his 1924 paper on combinators, there's only one other known piece of published output from Moses Schönfinkel: a paper coauthored with Paul Bernays "On the Decision Problem of Mathematical Logic":

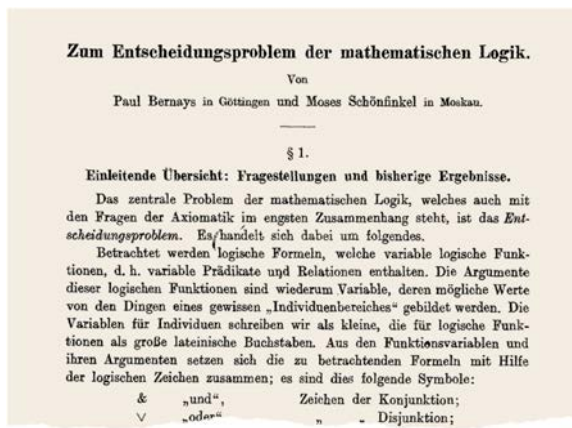

It's actually much more widely cited than Schönfinkel's 1924 combinator paper, but it's vastly less visionary and ultimately much less significant; it's really about a technical point in mathematical logic.

About halfway through the paper it has a note:

Form der Anzahl... ...g ... den Individuenbereich.

Die im folgenden mitgeteilten Überlegungen sind durch Vorlesungen von Hilbert über mathematische Logik angeregt worden und liegen schon um mehrere Jahre zurück. Die Durchführung des Entscheidungsverfahrens

350 P. Bernays und M. Schönfinkel.

im Falle einer einzigen auftretenden Funktion $F(x, y)$ rührt von M. Schönfinkel her, der zuerst das Problem in Angriff nahm[7], von P. Bernays die Ausdehnung der Methode auf mehrere logische Funktionen, sowie die Abfassung der vorliegenden Arbeit.

"The following thoughts were inspired by Hilbert's lectures on mathematical logic and date back several years. The decision procedure for a single function $F(x, y)$ was derived by M. Schönfinkel, who first tackled the problem; P. Bernays extended the method to several logical functions, and also wrote the current paper."

The paper was submitted on March 24, 1927. But in the records of the German Mathematicians Association we find a listing of another talk at the Göttingen Mathematical Society: December 6, 1921, P. Bernays and M. Schönfinkel, "Das *Entscheidungsproblem* im Logikkalkul". So the paper had a long gestation period, and (as the note in the paper suggests) it basically seems to have fallen to Bernays to get it written, quite likely with little or no communication with Schönfinkel.

So what else do we know about it? Well, remarkably enough, the Bernays archive contains two notebooks (the paper kind!) by Moses Schönfinkel that are basically an early draft of the paper (with the title already being the same as it finally was, but with Schönfinkel alone listed as the author):

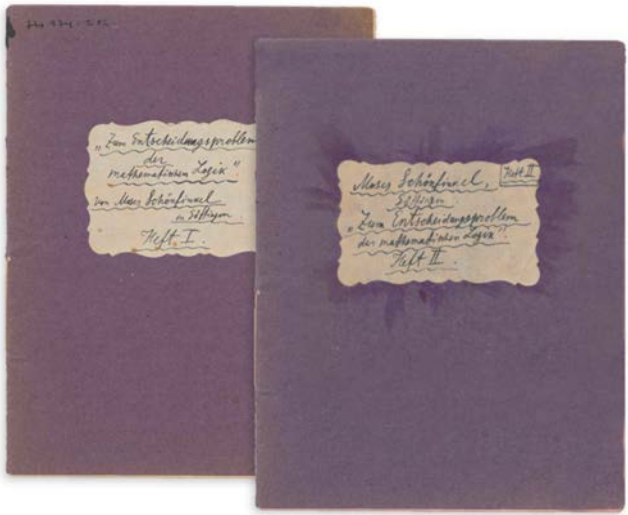

*ETH Zurich, Bernays Archive, Hs. 974: 282*

These notebooks are basically our best window into the front lines of Moses Schönfinkel's work. They aren't dated as such, but at the end of the second notebook there's a byline of sorts, that lists his street address in Göttingen—and we know he lived at that address from September 1922 until March 1924:

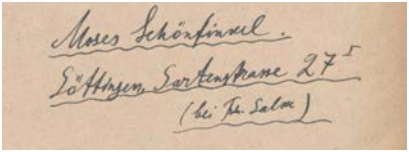

OK, so what's in the notebooks? The first page might indicate that the notebooks were originally intended for a different purpose. It's just a timetable of lectures:

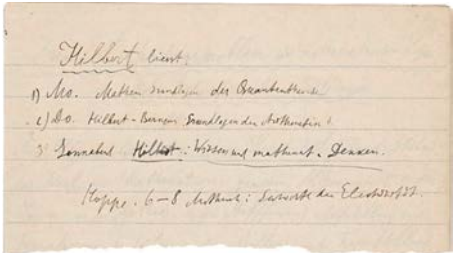

"Hilbert lectures: Monday: Mathematical foundations of quantum theory; Thursday: Hilbert–Bernays: Foundations of arithmetic; Saturday: Hilbert: Knowledge and mathematical thinking". (There's also a slightly unreadable note that seems to say "Hoppe. 6–8… electricity", perhaps referring to Edmund Hoppe, who taught physics in Göttingen, and wrote a history of electricity.)

But then we're into 15 pages (plus 6 in the other notebook) of content, written in essentially perfect German, but with lots of parentheticals of different possible word choices:

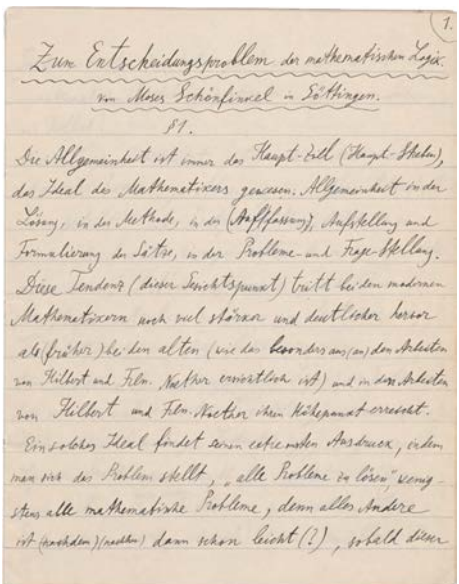

The final paper as coauthored with Bernays begins:

"The central problem of mathematical logic, which is also closely connected to its axiomatic foundations, is the decision problem [*Entscheidungsproblem*]. And it deals with the following. We have logical formulas which contain logic functions, predicates, ..."

Schönfinkel's version begins considerably more philosophically (here with a little editing for clarity):

"Generality has always been the main goal—the ideal of the mathematician. Generality in the solution, in the method, in the concept and formulation of the theorem, in the problem and question. This tendency is even more pronounced and clearer with modern mathematicians than with earlier ones, and reaches its high point in the work of Hilbert and Ms. Noether. Such an ideal finds its most extreme expression when one faces the problem of 'solving all problems'— at least all mathematical problems, because everything else after is easy, as soon as this 'Gordian Knot' is cut (because the world is written in 'mathematical letters' according to Hilbert).

In just the previous century mathematicians would have been extremely skeptical and even averse to such fantasies... But today's mathematician has already been trained and tested in the formal achievements of modern mathematics and Hilbert's axiomatics, and nowadays one has the courage and the boldness to dare to touch this question as well. We owe to mathematical logic the fact that we are able to have such a question at all.

From Leibniz's bold conjectures, the great logician-mathematicians went step by step in pursuit of this goal, in the systematic structure of mathematical logic: Boole (discoverer of the logical calculus), (Bolzano?), Ernst Schröder, Frege, Peano, Ms. Ladd-Franklin, the two Peirces, Sheffer, Whitehead, Couturat, Huntington, Padoa, Shatunovsky, Sleshinsky, Kagan, Poretsky, Löwenheim, Skolem, ... and their numerous students, collaborators and contemporaries ... until in 1910–1914 'the system' by Bertrand Russell and Whitehead appeared— the famous 'Principia Mathematica'—a mighty titanic work, a large system. Finally came our knowledge of logic from Hilbert's lectures on (the algebra of) logic (-calculus) and, following on from this, the groundbreaking work of Hilbert's students: Bernays and Behmann.

The investigations of all these scholars and researchers have led (in no uncertain terms) to the fact that it has become clear that actual mathematics represents a branch of logic. ... This emerges most clearly from the treatment and conception of mathematical logic that Hilbert has given. And now, thanks to Hilbert's approach, we can (satisfactorily) formulate the great decision problem of mathematical logic."

We learn quite a bit about Schönfinkel from this. Perhaps the most obvious thing is that he was a serious fan of Hilbert and his approach to mathematics (with a definite shout-out to "Ms. Noether"). It's also interesting that he refers to Bernays and Behmann as "students" of Hilbert. That's pretty much correct for Behmann. But Bernays (as we'll see soon) was more an assistant or colleague of Hilbert's than a student.

It gives interesting context to see Schönfinkel rattle off a sequence of contributors to what he saw as the modern view of mathematical logic. He begins—quite rightly I think—mentioning "Leibniz's bold conjectures". He's not sure whether Bernard Bolzano fits (and neither am I). Then he lists Schröder, Frege and Peano—all pretty standard choices, involved in building up the formal structure of mathematical logic.

Next he mentions Christine Ladd-Franklin. At least these days, she's not particularly well known, but she had been a mathematical logic student of Charles Peirce, and in 1881 she'd written a paper about the "Algebra of Logic" which included a truth table, a solid 40 years before Post or Wittgenstein. (In 1891 she had also worked in Göttingen on color vision with the experimental psychologist Georg Müller—who was still there in 1921.) It's notable that Schönfinkel mentions Ladd-Franklin ahead of the father-and-son Peirces. Next we see Sheffer, who Schönfinkel quotes in connection with Nand in his combinator paper. (No doubt unbeknownst to Schönfinkel, Henry Sheffer—who spent most of his life in the US— was also born in Ukraine ["near Odessa", his documents said], and was also Jewish, and was just 6 years older than Schönfinkel.) I'm guessing Schönfinkel mentions Whitehead next in connection with universal algebra, rather than his later collaboration with Russell.

Next comes Louis Couturat, who frankly wouldn't have made my list for mathematical logic, but was another "algebra of logic" person, as well as a Leibniz fan, and developer of the Ido language offshoot from Esperanto. Huntington was involved in the axiomatization of Boolean algebra; Padoa was connected to Peano's program. Shatunovsky, Sleshinsky and Kagan were all professors of Schönfinkel's in Odessa (as mentioned above), concerned in various ways with foundations of mathematics. Platon Poretsky I must say I had never heard of before; he seems to have done fairly technical work on propositional logic. And finally Schönfinkel lists Löwenheim and Skolem, both of whom are well known in mathematical logic today.

I consider it rather wonderful that Schönfinkel refers to Whitehead and Russell's *Principia Mathematica* as a "titanic work" (*Titanenwerk*). The showy and "overconfident" *Titanic* had come to grief on its iceberg in 1912, somehow reminiscent of *Principia Mathematica*, eventually coming to grief on Gödel's theorem.

At first it might just seem charming—particularly in view of his brother's comment that "[Moses] is helpless and impractical in this material world"—to see Schönfinkel talk about how after one's solved all mathematical problems, then solving *all problems* will be easy, explaining that, after all, Hilbert has said that "the world is written in 'mathematical letters'". He says that in the previous century mathematicians wouldn't have seriously considered "solving everything", but now, because of progress in mathematical logic, "one has the courage and the boldness to dare to touch this question".

It's very easy to see this as naive and unworldly—the writing of someone who knew only about mathematics. But though he didn't have the right way to express it, Schönfinkel was actually onto something, and something very big. He talks at the beginning of his piece about generality, and about how recent advances in mathematical logic embolden one to pursue it. And in a sense he was very right about this. Because mathematical logic—through work like his—is what led us to the modern conception of computation, which really is successful in "talking about everything". Of course, after Schönfinkel's time we learned about Gödel's theorem and computational irreducibility, which tell us that even though we may be able to talk about everything, we can never expect to "solve every problem" about everything.

But back to Schönfinkel's life and times. The remainder of Schönfinkel's notebooks give the technical details of his solution to a particular case of the decision problem. Bernays obviously worked through these, adding more examples as well as some generalization. And Bernays cut out Schönfinkel's philosophical introduction, no doubt on the (probably correct) assumption that it would seem too airy-fairy for the paper's intended technical audience.

So who was Paul Bernays? Here's a picture of him from 1928:

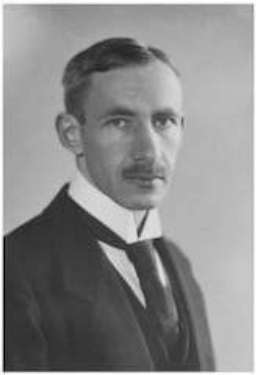

Bernays was almost exactly the same age as Schönfinkel (he was born on October 17, 1888—in London, where there was no calendar issue to worry about). He came from an international business family, was a Swiss citizen and grew up in Paris and Berlin. He studied math, physics and philosophy with a distinguished roster of professors in Berlin and Göttingen, getting his PhD in 1912 with a thesis on analytic number theory.

After his PhD he went to the University of Zurich, where he wrote a habilitation (on complex analysis), and became a privatdozent (yes, with the usual documentation, that can still be found), and an assistant to Ernst Zermelo (of ZFC set theory fame). But in 1917 Hilbert visited Zurich and soon recruited Bernays to return to Göttingen. In Göttingen, for apparently bureaucratic reasons, Bernays wrote a second habilitation, this time on the axiomatic structure of *Principia Mathematica* (again, all the documentation can still be found). Bernays was also hired to work as a "foundations of math assistant" to Hilbert. And it was presumably in that capacity that he—along with Moses Schönfinkel—wrote the notes for Hilbert's 1920 course on mathematical logic.

Unlike Schönfinkel, Bernays followed a fairly standard—and successful—academic track. He became a professor in Göttingen in 1922, staying there until he was dismissed (because of partially Jewish ancestry) in 1933—after which he moved back to Zurich, where he stayed and worked very productively, mostly in mathematical logic (von Neumann–Bernays–Godel set theory, etc.), until he died in 1977.

Back when he was in Göttingen one of the things Bernays did with Hilbert was to produce the two-volume classic *Grundlagen der Mathematik* (*Foundations of Mathematics*). So did the *Grundlagen* mention Schönfinkel? It has one mention of the Bernays–Schönfinkel paper, but no direct mention of combinators. However, there is one curious footnote:

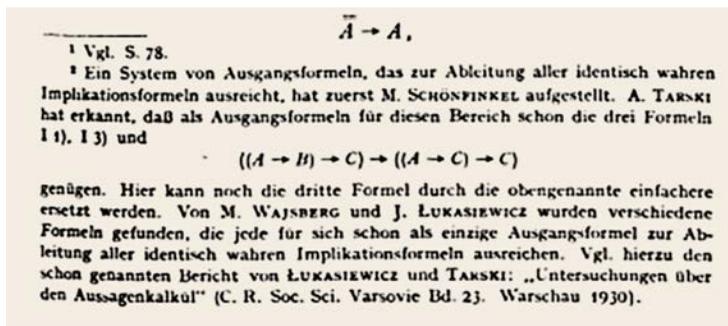

This starts "A system of axioms that is sufficient to derive all true implicational formulas was first set up by M. Schönfinkel…", then goes on to discuss work by Alfred Tarski. So do we have evidence of something else Schönfinkel worked on? Probably.

In ordinary logic, one starts from an axiom system that gives relations, say about A<small>ND</small>, O<small>R</small> and N<small>OT</small>. But, as Sheffer established in 1910, it's also possible to give an axiom system purely in terms of N<small>AND</small> (and, yes, I'm proud to say that I found the very simplest such axiom system in 2000). Well, it's also possible to use other bases for logic. And this footnote is about using I<small>MPLIES</small> as the basis. Actually, it's implicational calculus, which isn't as strong as ordinary logic, in the sense that it only lets you prove some of the theorems. But there's a question again: what are the possible axioms for implicational calculus?

Well, it seems that Schönfinkel found a possible set of such axioms, though we're not told what they were; only that Tarski later found a simpler set. (And, yes, I looked for the simpler axiom systems for implicational calculus in 2000, but didn't find any.) So again we see Schönfinkel in effect trying to explore the lowest-level foundations of mathematical logic, though we don't know any details.

So what other interactions did Bernays have with Schönfinkel? There seems to be no other information in Bernays's archives. But I have been able to get a tiny bit more information. In a strange chain of connections, someone who's worked on Mathematica and Wolfram Language since 1987 is Roman Maeder. And Roman's thesis advisor (at ETH Zurich) was Erwin Engeler—who was a student of Paul Bernays. Engeler (who is now in his 90s) worked for many years on combinators, so of course I had to ask him what Bernays might have told him about Schönfinkel. He told me he recalled only two conversations. He told me he had the impression that Bernays found Schönfinkel a difficult person. He also said he believed that the last time Bernays saw Schönfinkel it was in Berlin, and that Schönfinkel was somehow in difficult circumstances. Any such meeting in Berlin would have had to be before 1933. But try as we might to track it down, we haven't succeeded.

# To Moscow and Beyond…

In the space of three days in March 1924 Moses Schönfinkel—by then 35 years old—got his paper on combinators submitted to *Mathematische Annalen*, got a reference for himself sent out, and left for Moscow. But why did he go to Moscow? We simply don't know.

A few things are clear, though. First, it wasn't difficult to get to Moscow from Göttingen at that time; there was pretty much a direct train there. Second, Schönfinkel presumably had a valid Russian passport (and, one assumes, didn't have any difficulties from not having served in the Russian military during World War I).

One also knows that there was a fair amount of intellectual exchange and travel between Göttingen and Moscow. The very same volume of *Mathematische Annalen* in which Schönfinkel's paper was published has three (out of 19) authors in addition to Schönfinkel listed as being in Moscow: Pavel Alexandroff, Pavel Urysohn and Aleksandr Khintchine. Interestingly, all of these people were at Moscow State University.

And we know there was more exchange with that university. Nikolai Luzin, for example, got his PhD in Göttingen in 1915, and went on to be a leader in mathematics at Moscow State University (until he was effectively dismissed by Stalin in 1936). And we know that for example in 1930, Andrei Kolmogorov, having just graduated from Moscow State University, came to visit Hilbert.

Did Schönfinkel go to Moscow State University? We don't know (though we haven't yet been able to access any archives that may be there).

Did Schönfinkel go to Moscow because he was interested in communism? Again, we don't know. It's not uncommon to find mathematicians ideologically sympathetic to at least the theory of communism. But communism doesn't seem to have particularly been a thing in the mathematics or general university community in Göttingen. And indeed when Ludwig Gossmann was arrested in 1933, investigations of who he might have recruited into communism didn't find anything of substance.

Still, as I'll discuss later, there is a tenuous reason to think that Schönfinkel might have had some connection to Leon Trotsky's circle, so perhaps that had something to do with him going to Moscow—though it would have been a bad time to be involved with Trotsky, since by 1925 he was already out of favor with Stalin.

A final theory is that Schönfinkel might have had relatives in Moscow; at least it looks as if some of his Lurie cousins ended up there.

But realistically we don't know. And beyond the bylines on the journals, we don't really have any documentary evidence that Schönfinkel was in Moscow. However, there is one more data point, from November 1927 (8 months after the submission of Schönfinkel's paper with Bernays). Pavel Alexandroff was visiting Princeton University, and when Haskell Curry (who we'll meet later) asked him about Schönfinkel he was apparently told that "Schönfinkel has... gone insane and is now in a sanatorium & will probably not be able to work any more."

Ugh! What happened? Once again, we don't know. Schönfinkel doesn't seem to have ever been "in a sanatorium" while he was in Göttingen; after all, we have all his addresses, and none of them were sanatoria. Maybe there's a hint of something in Schönfinkel's brother's letter to Hilbert. But are we really sure that Schönfinkel actually suffered from mental illness? There's a bunch of hearsay that says he did. But then it's a common claim that logicians who do highly abstract work are prone to mental illness (and, well, yes, there are a disappointingly large number of historical examples).

Mental illness wasn't handled very well in the 1920s. Hilbert's only child, his son Franz (who was about five years younger than Schönfinkel), suffered from mental illness, and after a delusional episode that ended up with him in a clinic, David Hilbert simply said "From now on I have to consider myself as someone who does not have a son". In Moscow in the 1920s—despite some political rhetoric—conditions in psychiatric institutions were probably quite poor, and there was for example quite a bit of use of primitive shock therapy (though not yet electroshock). It's notable, by the way, that Curry reports that Alexandroff described Schönfinkel as being "in a sanatorium". But while at that time the word "sanatorium" was being used in the US as a better term for "insane asylum", in Russia it still had more the meaning of a place for a rest cure. So this still doesn't tell us if Schönfinkel was in fact "institutionalized"—or just "resting". (By the way, if there was mental illness involved, another connection for Schönfinkel that doesn't seem to have been made is that Paul Bernays's first cousin once removed was Martha Bernays, wife of Sigmund Freud.)

Whether or not he was mentally ill, what would it have been like for Schönfinkel in what was then the Soviet Union in the 1920s? One thing is that in the Soviet system, everyone was supposed to have a job. So Schönfinkel was presumably employed doing something—though we have no idea what. Schönfinkel had presumably been at least somewhat involved with the synagogue in Göttingen (which is how the rabbi there knew to tell his brother he was in bad shape). There was a large and growing Jewish population in Moscow in the 1920s, complete with things like Yiddish newspapers. But by the mid 1930s it was no longer so comfortable to be Jewish in Moscow, and Jewish cultural organizations were being shut down.

By the way, in the unlikely event that Schönfinkel was involved with Trotsky, there could have been trouble even by 1925, and certainly by 1929. And it's notable that it was a common tactic for Stalin (and others) to claim that their various opponents were "insane".

So what else do we know about Schönfinkel in Moscow? It's said that he died there in 1940 or 1942, aged 52–54. Conditions in Moscow wouldn't have been good then; the so-called Battle of Moscow occurred in the winter of 1941. And there are various stories told about Schönfinkel's situation at that time.

The closest to a primary source seems to be a summary of mathematical logic in the Soviet Union, written by Sofya Yanovskaya in 1948. Yanovskaya was born in 1896 (so 8 years after Schönfinkel), and grew up in Odessa. She attended the same university there as Schönfinkel, studying mathematics, though arrived five years after Schönfinkel graduated. She had many of the same professors as Schönfinkel, and, probably like Schönfinkel, was particularly influenced by Shatunovsky. When the Russian Revolution happened, Yanovskaya went "all in", becoming a serious party operative, but eventually began to teach, first at the Institute of Red Professors, and then from 1925 at Moscow State University—where she became a major figure in mathematical logic, and was eventually awarded the Order of Lenin.

One might perhaps have thought that mathematical logic would be pretty much immune to political issues. But the founders of communism had talked about mathematics, and there was a complex debate about the relationship between Marxist–Leninist ideology and formal ideas in mathematics, notably the Law of Excluded Middle. Sofya Yanovskaya was deeply involved, initially in trying to "bring mathematics to heel", but later in defending it as a discipline, as well as in editing Karl Marx's mathematical writings.

It's not clear to what extent her historical writings were censored or influenced by party considerations, but they certainly contain lots of good information, and in 1948 she wrote a paragraph about Schönfinkel:

закон исключенного третьего теряет силу.
7. Существенную роль в дальнейшем развитии математической логики сыграла работа М. И. Ш е й н ф и н к е л я [1]. Этот блестящий ученик С. О. Ш а т у н о в с к о г о, к сожалению, рано выбыл из строя. (Заболев душевно, М. И. Ш е й н ф и н к е л ь умер в Москве в 1942 г.). Работа, о которой идёт речь, была выполнена им в 1920 г., но опубликована только в 1924 г. в литературном оформлении Бемана. Непосредственной целью

"The work of M. I. Sheinfinkel played a substantial role in the further development of mathematical logic. This brilliant student of S. O. Shatunovsky, unfortunately, left us early. (After getting mentally ill [заболев душевно], M. I. Sheinfinkel passed away in Moscow in 1942.) He did the work mentioned here in 1920, but only published it in 1924, edited by Behmann."

Unless she was hiding things, this quote doesn't make it sound as if Yanovskaya knew much about Schönfinkel. (By the way, her own son was apparently severely mentally ill.) A student of Jean van Heijenoort (who we'll encounter later) named Irving Anellis did apparently in the 1990s ask a student of Yanovskaya's whether Yanovskaya had known Schönfinkel. Apparently he responded that unfortunately nobody had thought to ask her that question before she died in 1966.

What else do we know? Nothing substantial. The most extensively embellished story I've seen about Schönfinkel appears in an anonymous comment on the talk page for the Wikipedia entry about Schönfinkel:

"William Hatcher, while spending time in St Petersburg during the 1990s, was told by Soviet mathematicians that Schönfinkel died in wretched poverty, having no job and but one room in a collective apartment. After his death, the rough ordinary people who shared his apartment burned his manuscripts for fuel (WWII was raging). The few Soviet mathematicians around 1940 who had any discussions with Schönfinkel later said that those mss reinvented a great deal of 20th century mathematical logic. Schönfinkel had no way of accessing the work of Turing, Church, and Tarski, but had derived their results for himself. Stalin did not order Schönfinkel shot or deported to Siberia, but blame for Schönfinkel's death and inability to publish in his final years can be placed on Stalin's doorstep. 202.36.179.65 06:50, 25 February 2006 (UTC)"

William Hatcher was a mathematician and philosopher who wrote extensively about the Bahá'í Faith and did indeed spend time at the Steklov Institute of Mathematics in Saint Petersburg in the 1990s—and mentioned Schönfinkel's technical work in his writings. People I've asked at the Steklov Institute do remember Hatcher, but don't know anything about what it's claimed he was told about Schönfinkel. (Hatcher died in 2005, and I haven't been successful at getting any material from his archives.)

So are there any other leads? I did notice that the IP address that originated the Wikipedia comment is registered to the University of Canterbury in New Zealand. So I asked people there and in the New Zealand foundations of math scene. But despite a few "maybe so-and-so wrote that" ideas, nobody shed any light.

OK, so what about at least a death certificate for Schönfinkel? Well, there's some evidence that the registry office in Moscow has one. But they tell us that in Russia only direct relatives can access death certificates....

## Other Schönfinkels...

So far as we know, Moses Schönfinkel never married, and didn't have children. But he did have a brother, Nathan, who we encountered earlier in connection with the letter he wrote about Moses to David Hilbert. And in fact we know quite a bit about Nathan Scheinfinkel (as he normally styled himself). Here's a biographical summary from 1932:

Scheinfinkel, Nathan, Dr. med., P.-D., Neufeldstr. 5a, Bern (geb. 13. IX. 93 in Ekaterinoslaw, Ukraine). Nat.: Ausländer. — Stud. Gymn. Ekaterinoslaw, Univ. Bern. 20 med. Dr.-Examen; seit 22 Assistent a. Berner Physiolog. Institut (Ausbildung unt. Prof. Asher nach d. biochem. u. biophysikal. Richtung im Geiste d. modernen Forschung.) Seit 29 P.-D. f. Physiologie a. d. Univ. Bern. — V.: Sämtl. Arbeiten beschäftigen sich mit den Problemen: Ermüdung, Sauerstoffmangel, Herzarbeit u. Wirkungsweise d. Herznerven (ersch. i. d. Zeitschr. f. Biologie).

*Deutsches Biographisches Archiv, II 1137, 103*

The basic story is that he was about five years younger than Moses, and went to study medicine at the University of Bern in Switzerland in April 1914 (i.e. just before World War I began). He got his MD in 1920, then got his PhD on "Gas Exchange and Metamorphosis of Amphibian Larvae after Feeding on the Thyroid Gland or Substances Containing Iodine" in 1922. He did subsequent research on the electrochemistry of the nervous system, and in 1929 became a privatdozent—with official "license to teach" documentation:

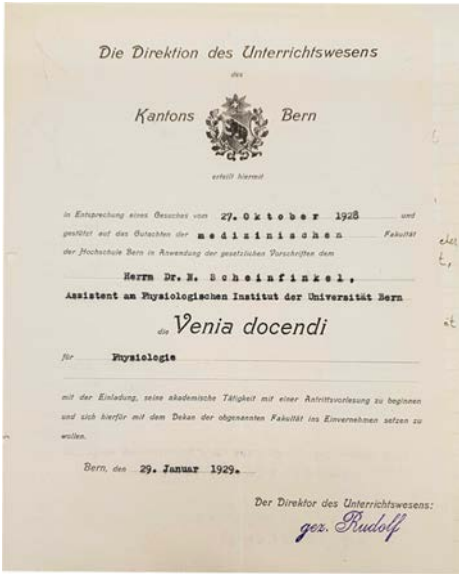

*Staatsarchiv des Kantons Bern, BB IIIb 557 Scheinfinkel N.*

(In a piece of bizarre small-worldness, my grandfather, Max Wolfram, also got a PhD in the physiology [veterinary medicine] department at the University of Bern [studying the function of the thymus gland], though that was in 1909, and presumably he had left before Nathan Scheinfinkel arrived.)

But in any case, Nathan Scheinfinkel stayed at Bern, eventually becoming a professor, and publishing extensively, including in English. He became a Swiss citizen in 1932, with the official notice stating:

> "Scheinfinkel, Nathan. Son of Ilia Gerschow and Mascha [born] Lurie, born in Yekaterinoslav, Russia, September 13, 1893 (old style). Doctor of medicine, residing in Bern, Neufeldstrasse 5a, husband of Raissa [born] Neuburger."

In 1947, however, he moved to become a founding professor in a new medical school in Ankara, Turkey. (Note that Turkey, like Switzerland, had been neutral in World War II.) In 1958 he moved again, this time to found the Institute of Physiology at Ege University in Izmir, Turkey, and then at age 67, in 1961, he retired and returned to Switzerland.

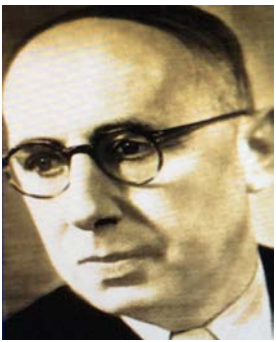

Did Nathan Scheinfinkel have children (whose descendents, at least, might know something about "Uncle Moses")? It doesn't seem so. We tracked down Nuran Harirî, now an emeritus professor, but in the 1950s a young physiology resident at Ege University responsible for translating Nathan Scheinfinkel's lectures into Turkish. She said that Nathan Scheinfinkel was at that point living in campus housing with his wife, but she never heard mention of any children, or indeed of any other family members.

What about any other siblings? Amazingly, looking through handwritten birth records from Ekaterinoslav, we found one! Debora Schönfinkel, born December 22, 1889 (i.e. January 3, 1890, in the modern calendar):

So Moses Schönfinkel had a younger sister, as well as a younger brother. And we even know that his sister graduated from high school in June 1907. But we don't know anything else about her, or about other siblings. We know that Schönfinkel's mother died in 1936, at the age of 74.

Might there have been other Schönfinkel relatives in Ekaterinoslav? Perhaps, but it's unlikely they survived World War II—because in one of those shocking and tragic pieces of history, over a four-day period in February 1942 almost the whole Jewish population of 30,000 was killed.

Could there be other Schönfinkels elsewhere? The name is not common, but it does show up (with various spellings and transliterations), both before and after Moses Schönfinkel. There's a Scheinfinkel Russian revolutionary buried in the Kremlin Wall; there was a Lovers of Zion delegate Scheinfinkel from Ekaterinoslav. There was a Benjamin Scheinfinkel in New York City in the 1940s; a Shlomo Scheinfinkel in Haifa in the 1930s. There was even a certain curiously named Bas Saul Haskell Scheinfinkel born in 1875. But despite quite a bit of effort, I've been unable to locate any living relative of Moses Schönfinkel. At least so far.

## Haskell Curry

What happened with combinators after Schönfinkel published his 1924 paper? Initially, so far as one can tell, nothing. That is, until Haskell Curry found Schönfinkel's paper in the library at Princeton University in November 1927—and launched into a lifetime of work on combinators.

Who was Haskell Curry? And why did he know to care about Schönfinkel's paper?

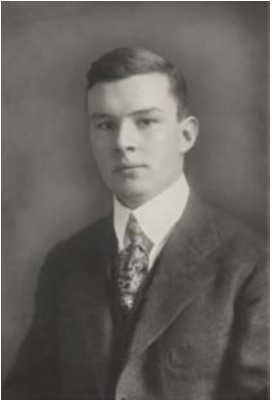

Haskell Brooks Curry was born on September 12, 1900, in a small town near Boston, MA. His parents were both elocution educators, who by the time Haskell Curry was born were running the School of Expression (which had evolved from his mother's Boston-based School of Elocution and Expression). (Many years later, the School of Expression would evolve into Curry College in Waltham, Massachusetts—which happens to be where for several years we held our Wolfram Summer School, often noting the "coincidence" of names when combinators came up.)

Haskell Curry went to college at Harvard, graduating in mathematics in 1920. After a couple of years doing electrical engineering, he went back to Harvard, initially working with Percy Bridgman, who was primarily an experimental physicist, but was writing a philosophy of science book entitled *The Logic of Modern Physics*. And perhaps through this Curry got introduced to Whitehead and Russell's *Principia Mathematica*.

But in any case, there's a note in his archive about *Principia Mathematica* dated May 20, 1922:

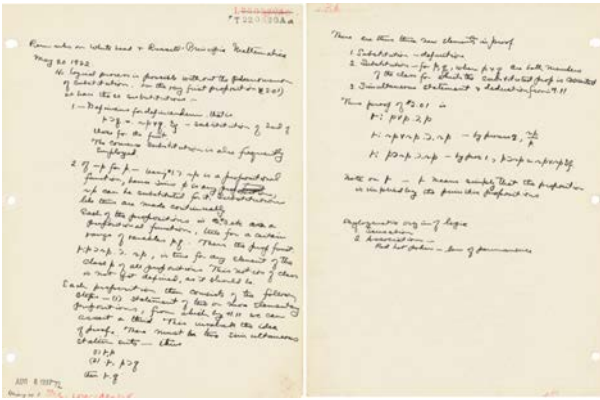

Curry seems—perhaps like an electrical engineer or a "pre-programmer"—to have been very interested in the actual process of mathematical logic, starting his notes with: "No logical process is possible without the phenomenon of substitution." He continued, trying to break down the process of substitution.

But then his notes end, more philosophically, and perhaps with "expression" influence: "Phylogenetic origin of logic: 1. Sensation; 2. Association: Red hot poker–law of permanence".

At Harvard Curry started working with George Birkhoff towards a PhD on differential equations. But by 1927–8 he had decided to switch to logic, and was spending a year as an instructor at Princeton. And it was there—in November 1927—that he found Schönfinkel's paper. Preserved in his archives are the notes he made:

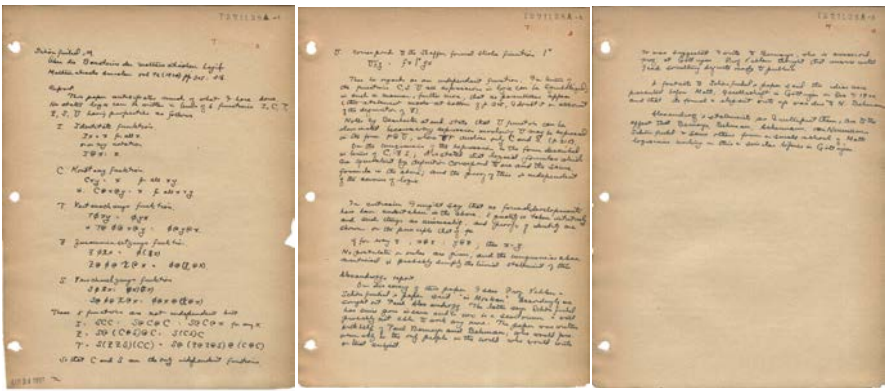

At the top there's a date stamp of November 28, 1927. Then Curry writes: "This paper anticipates much of what I have done"—then launches into a formal summary of Schönfinkel's paper (charmingly using *f@x* to indicate function application—just as we do in Wolfram Language, except his is left associative...).

He ends his "report" with "In criticism I might say that no formal development have been undertaken in the above. Equality is taken intuitively and such things as universality, and proofs of identity are shown on the principle that if for every $z$, $x@z : y@z$ then $x=y$ ...."

But then there's another piece:

"On discovery of this paper I saw Prof. Veblen. Schönfinkel's paper said 'in Moskau'. Accordingly we sought out Paul Alexandroff. The latter says Schönfinkel has since gone insane and is now in a sanatorium & will probably not be able to work any more. The paper was written with help of Paul Bernays and Behman [sic]; who would presumably be the only people in the world who would write on that subject."

What was the backstory to this? Oswald Veblen was a math professor at Princeton who had worked on the axiomatization of geometry and was by then working on topology. Pavel Alexandroff (who we encountered earlier) was visiting from Moscow State University for the year, working on topology with Hopf, Lefschetz, Veblen and Alexander. I'm not quite sure why Curry thought Bernays and Behmann "would be the only people in the world who would write on that subject"; I don't see how he could have known.

Curry continues: "It was suggested I write to Bernays, who is *außerord.* prof. [long-term lecturer] at Göttingen." But then he adds—in depressingly familiar academic form: "Prof. Veblen thought it unwise until I had something definite ready to publish."

"A footnote to Schönfinkel's paper said the ideas were presented before Math Gesellschaft in Göttingen on Dec. 7, 1920 and that its formal and elegant [sic] write up was due to H. Behman". "Elegant" is a peculiar translation of "*stilistische*" that probably gives Behmann too much credit; a more obvious translation might be "stylistic".

Curry continues: "Alexandroff's statements, as I interpret them, are to the effect that Bernays, Behman, Ackermann, von Neumann, Schönfinkel & some others form a small school of math logicians working on this & similar topics in Göttingen."

And so it was that Curry resolved to study in Göttingen, and do his PhD in logic there. But before he left for Göttingen, Curry wrote a paper (published in 1929):

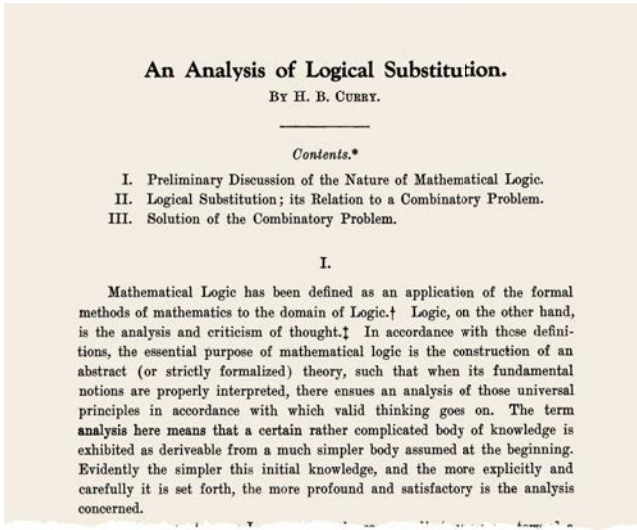

Already there's something interesting in the table of contents: the use of the word "combinatory", which, yes, in Curry's care is going to turn into "combinator".

The paper starts off reading a bit like a student essay, and one's not encouraged by a footnote a few pages in:

> "In writing the foregoing account I have naturally made use of any ideas I may have gleaned from reading the literature. The writings of Hilbert are fundamental in this connection. I hope that I have added clearness to certain points where the existing treatments are obscure." ["Clearness" not "clarity"?]

Then, towards the end of the "Preliminary Discussion" is this:

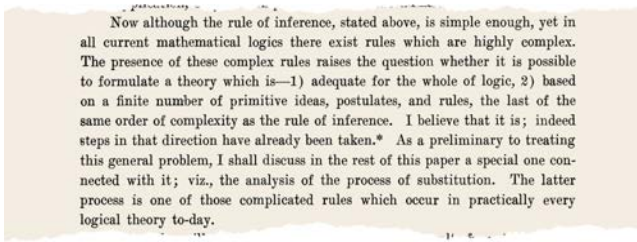

And the footnote says: "See the paper of Schönfinkel cited below". It's (so far as I know) the first-ever citation to Schönfinkel's paper!

On the next page Curry starts to give details. Curry starts talking about substitution, then says (in an echo of modern symbolic language design) this relates to the idea of "transform-ation of functions":

> A notion closely related to substitution is that of transformation of func-tions. Suppose we regard a function as having inherent in its definition a certain order of its variables. Then permuting these variables in any way, or making two or more of them alike, will produce new functions related to the old; let us call them transforms of the original function, and the opera-tions by which they are produced transformations. If we number the varia-bles consecutively 1, 2, 3, $\cdots$, then the transforms for a function of two variables will be—
>
> $$\phi(1, 2),\ \phi(2, 1),\ \phi(1, 1).$$
>
> For three variables there will be 13 transforms, for four variables 75, for five variables 541, etc. It is clear that the process of substituting a series of con-stants in an arbitrary manner (such that the total number of entities counting repetitions is $n$) into the original function is equivalent to the substitution of the same entities in a prescribed manner (viz., the first entity into the place of the first variable, the second into that of the second, etc.) into one of the transforms. The study of substitution is thus to a certain degree equivalent to the study of these transformations.
>
> An important step toward the analysis of this situation was made by M. Schönfinkel.* Starting, apparently, from the fact that every logical for-mula is a combination of constants—the variables being only apparent—he shows that neither the notions of propositional function (of various orders)
>
> ——————
> * "Ueber die Bausteine der mathematischen Logik," *Mathematische Annalen*, Vol. 92 (1924), pp. 305-316.

At first he's off talking about all the various combinatorial arrangements of variables, etc. But then he introduces Schönfinkel—and starts trying to explain in a formal way what Schönfinkel did. And even though he says he's talking about what one assumes is structural substitution, he seems very concerned about what equality means, and how Schönfinkel didn't quite define that. (And, of course, in the end, with universal computation, undecidabil-ity, etc. we know that the definition of equality wasn't really accessible in the 1920s.)

By the next page, here we are, *S* and *K* (Curry renamed Schönfinkel's *C*):

> CURRY: *An Analysis of Logical Substitution.* 371
>
> **III. POSTULATES.**
>    None.
>
> **IV. RULES.**
>
>    0. If $x$ and $y$ are entities, then $(xy)$ shall be an entity.
>
>    1. $(=\!=)$ shall have the properties of identity. These properties may be specified by a few simple rules; but in this treatment we shall not go into that detail. We shall treat $(=\!=)$ as if it were precisely the intuitive relation of equality.
>
>    2. If $x$ and $y$ are any entities, then
>
>    $$Kxy = x$$
>
>    3. If $x$, $y$, $z$ are entities, then
>
>    $$Sxyz = xz(yz)$$
>
>    4. If $X$ and $Y$ are combinations of $S$ and $K$, and if there exists an integer $n$ such that by application of the preceding rules we can formally reduce the expressions $Xx_1\ x_2\ \cdots\ x_n$ and $Yx_1\ x_2\ \cdots\ x_n$ to combinations of $x_1\ x_2\ \cdots\ x_n$ which have the same structure, then $X = Y$.
>
>    If the above primitive frame were a part of a general theory of logic, the term entity would include not only the various combinations of $S$ and $K$, but all the notions of logic as well. In the sequel we shall accordingly speak of the application of combinations $S$ and $K$ to various logical notions, and of the resulting notions to each other, just as if these notions had been adjoined to the above frame.
>
>    The *raison d'être* of the theory based on this frame is the following fact: Let $x_1$, $x_2$, $\cdots$, $x_n$ be any $n$ entities, and $X$ any combination of them con-At first he's imagining that the combinators have to be applied to something (i.e. $f[x]$ not just $f$). But by the next page he comes around to what Schönfinkel was doing in looking at "pure combinators":

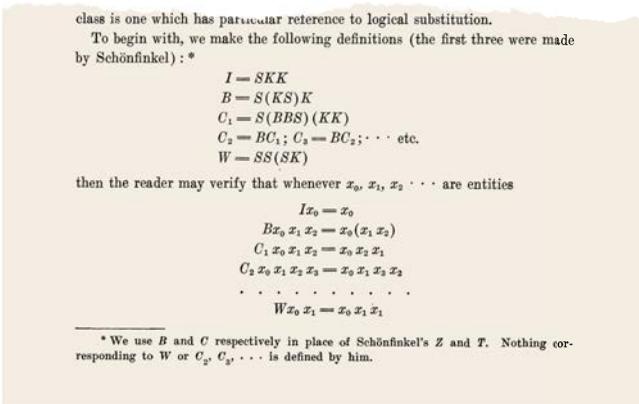

The rest of the paper is basically concerned with setting up combinators that can successively represent permutations—and it certainly would have been much easier if Curry had had a computer (and one could imagine minimal "combinator sorters" like minimal sorting networks):

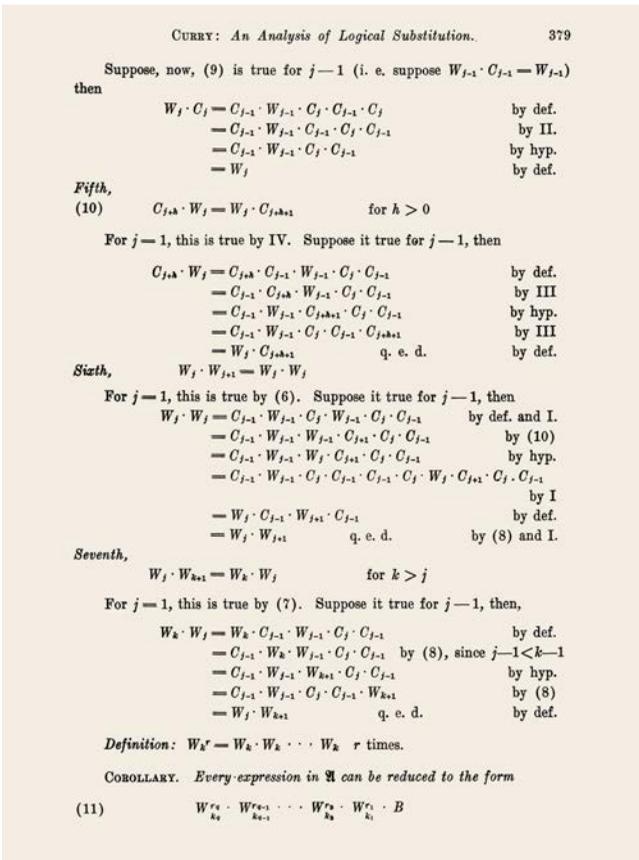

After writing this paper, Curry went to Göttingen—where he worked with Bernays. I must say that I'm curious what Bernays said to Curry about Schönfinkel (was it more than to Erwin Engeler?), and whether other people around Göttingen even remembered Schönfinkel, who by then had been gone for more than four years. In 1928, travel in Europe was open enough that Curry should have had no trouble going, for example, to Moscow, but there's no evidence he made any effort to reach out to Schönfinkel. But in any case, in Göttingen he worked on combinators, and over the course of a year produced his first official paper on "combinatory logic":

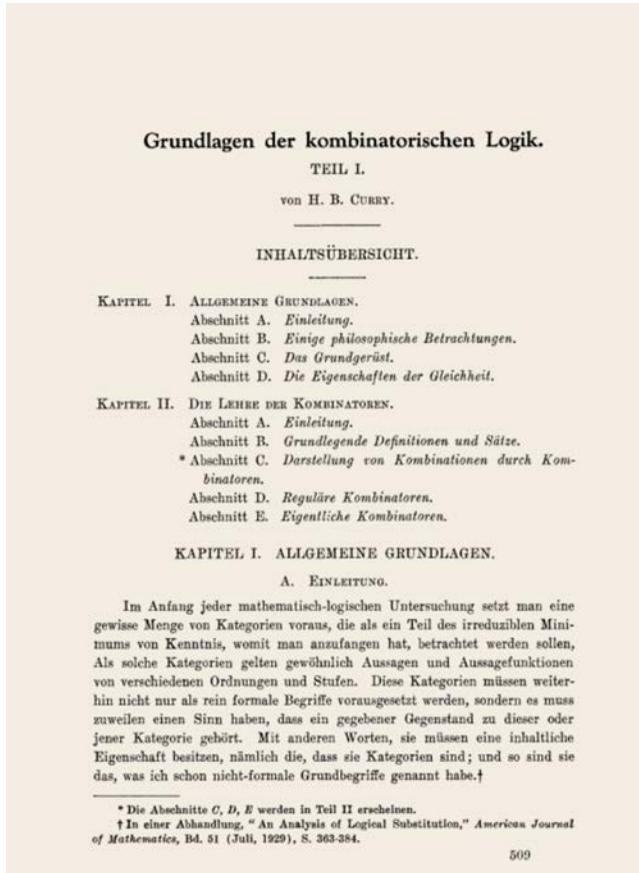

Strangely, the paper was published in an American journal—as the only paper not in English in that volume. The paper is more straightforward, and in many ways more "Schönfinkel like". But it was just the first of many papers that Curry wrote about combinators over the course of nearly 50 years.

Curry was particularly concerned with the "mathematicization" of combinators, finding and fixing problems with axioms invented for them, connecting to other formalisms (notably Church's lambda calculus), and generally trying to prove theorems about what combinators do. But more than that, Curry spread the word about combinators far and wide. And before long most people viewed him as "Mr. Combinator", with Schönfinkel at most a footnote.

In 1958, when Haskell Curry and Robert Feys wrote their book on *Combinatory Logic,* there's a historical footnote—that gives the impression that Curry "almost" had Schönfinkel's ideas before he saw Schönfinkel's paper in 1927:

**1. Historical statement**

The material in this chapter came mostly from [GKL], [PKR], and Rosser [MLV]. These have been supplemented from notes not previously published, and from specific suggestions made by various persons, notably Church and Bernays.

The combinators **B, C, I, K,** and **S** were introduced by Schönfinkel [BML] (see § 0D). This paper is still available for the discussion of the intuitive meaning of these combinators and the motivation for introducing them. (There is an error in the supplementary statement by Behmann at the end of the paper, in that it is not always possible to remove parentheses by **B** alone; Behmann has written that this error was pointed out by Boscovitch at an early date.)

The earliest work of Curry (till the fall of **1927**), which was done without knowledge of the work of Schönfinkel, used **B, C, W,** and **I** as primitive combinators. (These symbols were suggested by English, rather than German, names; '**B**' by 'substitution', since '**S**' might be used for several other purposes such as 'sum', 'successor', etc.; and '**W**' by a natural association of this letter with repetition.) When Schönfinkel was discovered in a literature search, **K** was added to the theory at once;

but **S** was regarded as a mere technicality until the development axiomatic theories in the 1940's (see § **6S1**).

The combinators **Bᵐ, Cₙ, Kₙ, Wₙ** were introduced in [GKL]; ₫ The combinators **Zₙ** were suggested, at least in principle, by Ch The connection between combinators and arithmetic was elabora and Rosser (cf. § **0C**); Rosser's [MLV] was especially significant i connection. The connections of **B** with products and powers was o

I have to say that I don't think that's a correct impression. What Schönfinkel did was much more singular than that. It's plausible to think that others (and particularly Curry) could have had the idea that there could be a way to go "below the operations of mathematical logic" and find more fundamental building blocks based on understanding things like the process of substitution. But the actuality of how Schönfinkel did it is something quite different—and something quite unique.

And when one sees Schönfinkel's *S* combinator: what mind could have come up with such a thing? Even Curry says he didn't really understand the significance of the *S* combinator until the 1940s.

I suppose if one's just thinking of combinatory logic as a formal system with a certain general structure then it might not seem to matter that things as simple as *S* and *K* can be the ultimate building blocks. But the whole point of what Schönfinkel was trying to do (as the title of his paper says) was to find the "building blocks of logic". And the fact that he was able to do it— especially in terms of things as simple as *S* and *K*—was a great and unique achievement. And not something that (despite all the good he did for combinators) Curry did.

## Schönfinkel Rediscovered

In the decade or so after Schönfinkel's paper appeared, Curry occasionally referenced it, as did Church and a few other closely connected people. But soon Schönfinkel's paper— and Schönfinkel himself—disappeared completely from view, and standard databases list no citations.

But in 1967 Schönfinkel's paper was seen again—now even translated into English. The venue was a book called *From Frege to Gödel: A Source Book in Mathematical Logic, 1879–1931.*

And there, sandwiched between von Neumann on transfinite numbers and Hilbert on "the infinite", is Schönfinkel's paper, in English, with a couple of pages of introduction by Willard Van Orman Quine. (And indeed it was from this book that I myself first became aware of Schönfinkel and his work.)

But how did Schönfinkel's paper get into the book? And do we learn anything about Schönfinkel from its appearance there? Maybe. The person who put the book together was a certain Jean van Heijenoort, who himself had a colorful history. Born in 1912, he grew up mostly in France, and went to college to study mathematics—but soon became obsessed with communism, and in 1932 left to spend what ended up being nearly ten years working as a kind of combination PR person and bodyguard for Leon Trotsky, initially in Turkey but eventually in Mexico. Having married an American, van Heijenoort moved to New York City, eventually enrolling in a math PhD program, and becoming a professor doing mathematical logic (though with some colorful papers along the way, with titles like "The Algebra of Revolution").

Why is this relevant? Well, the question is: how did van Heijenoort know about Schönfinkel? Perhaps it was just through careful scholarship. But just maybe it was through Trotsky. There's no real evidence, although it is known that during his time in Mexico, Trotsky did request a copy of *Principia Mathematica* (or was it his "PR person"?). But at least if there was a Trotsky connection it could help explain Schönfinkel's strange move to Moscow. But in the end we just don't know.

## What Should We Make of Schönfinkel?

When one reads about the history of science, there's a great tendency to get the impression that big ideas come suddenly to people. But my historical research—and my personal experience—suggest that that's essentially never what happens. Instead, there's usually a period of many years in which some methodology or conceptual framework gradually develops, and only then can the great idea emerge.

So with Schönfinkel it's extremely frustrating that we just can't see that long period of development. The records we have just tell us that Schönfinkel announced combinators on December 7, 1920. But how long had he been working towards them? We just don't know.

On the face of it, his paper seems simple—the kind of thing that could have been dashed off in a few weeks. But I think it's much more likely that it was the result of a decade of development—of which, through foibles of history, we now have no trace.

Yes, what Schönfinkel finally came up with is simple to explain. But to get to it, he had to cut through a whole thicket of technicality—and see the essence of what lay beneath. My life as a computational language designer has often involved doing very much this same kind of thing. And at the end of it, what you come up with may seem in retrospect "obvious". But to get there often requires a lot of hard intellectual work.

And in a sense what Schönfinkel did was the most impressive possible version of this. There were no computers. There was no ambient knowledge of computation as a concept. Yet Schönfinkel managed to come up with a system that captures the core of those ideas. And while he didn't quite have the language to describe it, I think he did have a sense of what he was doing—and the significance it could have.

What was the personal environment in which Schönfinkel did all this? We just don't know. We know he was in Göttingen. We don't think he was involved in any particularly official way with the university. Most likely he was just someone who was "around". Clearly he had some interaction with people like Hilbert and Bernays. But we don't know how much. And we don't really know if they ever thought they understood what Schönfinkel was doing.

Even when Curry picked up the idea of combinators—and did so much with it—I don't think he really saw the essence of what Schönfinkel was trying to do. Combinators and Schönfinkel are a strange episode in intellectual history. A seed sown far ahead of its time by a person who left surprisingly few traces, and about whom we know personally so little.

But much as combinators represent a way of getting at the essence of computation, perhaps in combinators we have the essence of Moses Schönfinkel: years of a life compressed to two "signs" (as he would call them) $S$ and $K$. And maybe if the operation we now call currying needs a symbol we should be using the "sha" character Ш from the beginning of Schönfinkel's name to remind us of a person about whom we know so little, but who planted a seed that gave us so much.

## Thanks

Many people and organizations have helped in doing research and providing material for this piece. Thanks particularly to Hatem Elshatlawy (fieldwork in Göttingen, etc.), Erwin Engeler (first-person history), Unal Goktas (Turkish material), Vitaliy Kaurov (locating Ukraine + Russia material), Anna & Oleg Marichev (interpreting old Russian handwriting), Nik Murzin (fieldwork in Moscow), Eila Stiegler (German translations), Michael Trott (interpreting German). Thanks also for input from Henk Barendregt, Semih Baskan, Metin Baştuğ, Cem Bozşahin, Jason Cawley, Jack Copeland, Nuran Harirî, Ersin Koylu, Alexander Kuzichev, Yuri Matiyasevich, Roman Maeder, Volker Peckhaus, Jonathan Seldin, Vladimir Shalack, Matthew Szudzik, Christian Thiel, Richard Zach. Particular thanks to the following archives and staff: Berlin State Library [Gabriele Kaiser], Bern University Archive [Niklaus Bütikofer], ETHZ (Bernays) Archive [Flavia Lanini, Johannes Wahl], Göttingen City Archive [Lena Uffelmann], Göttingen University [Katarzyna Chmielewska, Bärbel Mund, Petra Vintrová, Dietlind Willer].

# A Little Closer to Finding What Became of Moses Schönfinkel, Inventor of Combinators

For most big ideas in recorded intellectual history one can answer the question: "What became of the person who originated it?" But late last year I tried to answer that for Moses Schönfinkel, who sowed a seed for what's probably the single biggest idea of the past century: abstract computation and its universality.

I managed to find out quite a lot about Moses Schönfinkel. But I couldn't figure out what became of him. Still, I kept on digging. And it turns out I was able to find out more. So here's an update....

To recap a bit: Moses Schönfinkel was born in 1888 in Ekaterinoslav (now Dnipro) in what's now Ukraine. He went to college in Odessa, and then in 1914 went to Göttingen to work with David Hilbert. He didn't publish anything, but on December 7, 1920—at the age of 32—he gave a lecture entitled "Elemente der Logik" ("Elements of Logic") that introduced what are now called combinators, the first complete formalism for what we'd now call abstract computation. Then on March 18, 1924, with a paper based on his lecture just submitted for publication, he left for Moscow. And basically vanished.

It's said that he had mental health issues, and that he died in poverty in Moscow in 1940 or 1942. But we have no concrete evidence for either of these claims.

When I was researching this last year, I found out that Moses Schönfinkel had a younger brother Nathan Scheinfinkel (yes, he used a different transliteration of the Russian Шейнфинкель) who became a physiology professor at Bern in Switzerland, and later in Turkey. Late in the process, I also found out that Moses Schönfinkel had a younger sister Debora, who we could tell graduated from high school in 1907.

Moses Schönfinkel came from a Jewish merchant family, and his mother came from a quite prominent family. I suspected that there might be other siblings (Moses's mother came from a family of 8). And the first "new find" was that, yes, there were indeed two additional younger brothers. Here are the recordings of their births now to be found in the State Archives of the Dnipropetrovsk (i.e. Ekaterinoslav) Region:

So the complete complement of Шейнфинкель/Schönfinkel/Scheinfinkel children was (including birth dates both in their original Julian calendar form, and in their modern Gregorian form, and graduation dates in modern form):

| | *born (original calendar)* | *born (modern calendar)* | *graduated high school* |
|---|---|---|---|
| Moses (Моисей) | September 17, 1888 | September 29, 1888 | June 26, 1906 |
| Debora (Дебора) | December 22, 1889 | January 3, 1890 | June 25, 1907 |
| Nathan (Натан) | September 13, 1893 | September 25, 1893 | June 16, 1913 |
| Israel (Израиль) | December 5, 1894 | December 17, 1894 | June 16, 1913 |
| Gregory (Григорий) | April 30, 1899 | May 12, 1899 | ? 1917 |

And having failed to find out more about Moses Schönfinkel directly, plan B was to investigate his siblings.

I had already found out a fair amount about Nathan. He was married, and lived at least well into the 1960s, eventually returning to Switzerland. And most likely he had no children.

Debora we could find no trace of after her high-school graduation (we looked for marriage records, but they're not readily available for what we assume is the relevant time period).

By the way, rather surprisingly, we found nice (alphabetically ordered), printed class lists from the high-school graduations (apparently these were distributed to higher-education institutions across the Russian Empire so anyone could verify "graduation status", and were deposited in the archives of the education district, where they've now remained for more than a century):

(We can't find any particular trace of the 36 other students in the same group as Moses.)

OK, so what about the "newly found siblings", Israel and Gregory? Well, here we had a bit more luck.

For Israel we found these somewhat strange traces:

They are World War I hospital admission records from January and December 1916. Apparently Israel was a private in the 2nd Finnish Regiment (which—despite its name—by then didn't have any Finns in it, and in 1916 was part of the Russian 7th Army pushing west in southern Ukraine in the effort to retake Galicia). And the documents we have show that twice he ended up in a hospital in Pavlohrad (only about 40 miles from Ekaterinoslav, though in the opposite direction from where the 7th Army was) with some kind of (presumably not life-threatening) hernia-like problem.

But unfortunately, that's it. No more trace of Israel.

OK, what about the "baby brother", Gregory, 11 years younger than Moses? Well, he shows up in World War II records. We found four documents:

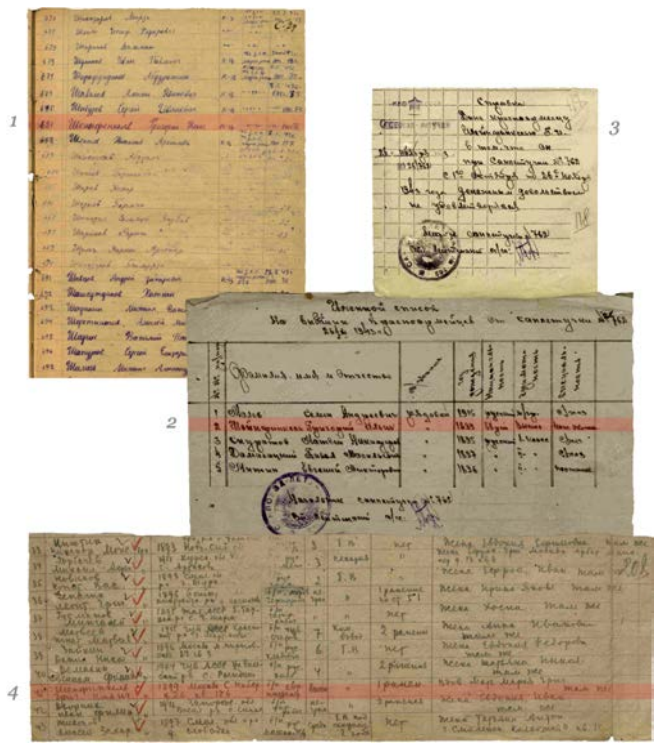

[Document #4](#) contains something interesting: an address for Gregory in 1944—in Moscow. Remember that Moses went to Moscow in 1924. And one of my speculations was that this was the result of some family connection there. Well, at least 20 years later (and probably also much earlier, as we'll see), his brother Gregory was in Moscow. So perhaps that's why Moses went there in 1924.

OK, but what story do these World War II documents tell about Gregory? Document #1 tells us that on July 27, 1943, Gregory arrived at the military unit designated 15 зсп 44 зсбр (15 ZSP 44 ZSBR) at transit point (i.e. basically "military address") 215 азсп 61А (215 AZSP 61A). It also tells us that he had the rank of private in the Red Army.

Sometime soon thereafter he was transferred to unit 206 ZSP. But unfortunately he didn't last long in the field. Around October 1, 1943, he was wounded (later, we learn he has "one wound"), and—as document #2 tells us—he was one of 5 people picked up by hospital train #762 (at transit point 206 зсп ЗапФ). On November 26, 1943, document #3 records that he was discharged from the hospital train (specifically, the document explains that he's not getting paid for the time he was on the hospital train). And, finally, document #4 records that on February 18, 1944—presumably after a period of assessment of his condition—he's discharged from the military altogether, returning to an address in Moscow.

OK, so first some military points. When Gregory arrived in the army in July 1943 he was assigned (as a reserve or "replacement") to the 44th Rifle Brigade (44 зсбр) in the 15th Rifle Division (15 зсп) in the 61st Army (61A)—presumably as part of reinforcements brought in after some heavy Soviet losses. Later he was transferred to the 206th Rifle Division in the 47th Army, which is where he was when he was wounded around October 1, 1943.

What was the general military situation then? In the summer of 1943 the major story was that the Soviets were trying to push the Germans back west, with the front pretty much along the Dnieper River in Ukraine—which, curiously enough, flows right through the middle of Ekaterinoslav. On October 4, 1943, here's how the *New York Times* presented things:

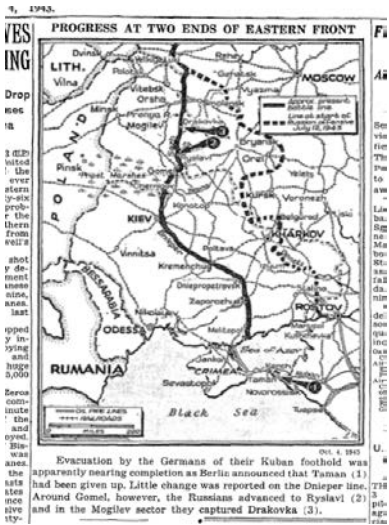

But military history being what it is, there's much more detailed information available. Here's a modern map showing troop movements involving the 47th Army in late September 1943:

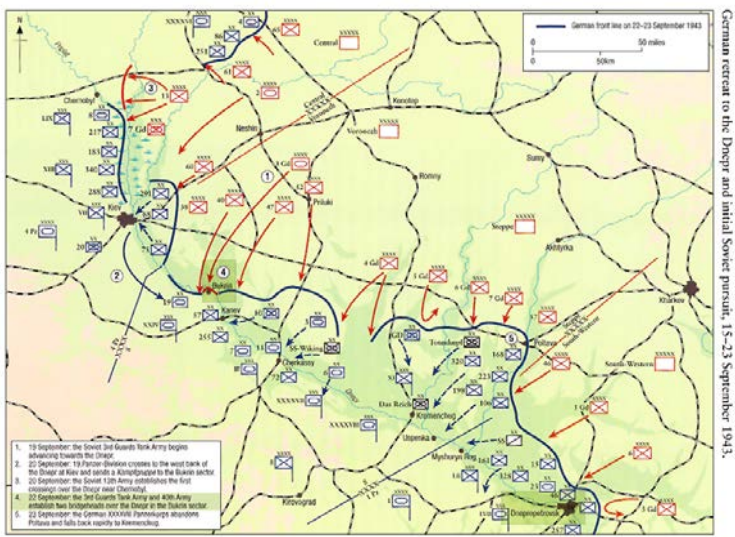

The Soviets managed to get more than 100,000 men across the Dnieper River, but there was intense fighting, and at the end of September the 206th Rifle Division (as part of the 47th Army) was probably involved in the later stages of the fight for the Bukrin Bridgehead. And this is probably where Gregory Schönfinkel was wounded.

After being wounded, he seems to have been taken to some kind of service area for the 206th Rifle Division (206 зсп ЗапФ), from which he was picked up by a hospital train (and, yes, it was actually a moving hospital, with lots of cars with red crosses painted on top).

But more significant in our quest for the story of Gregory Schönfinkel is other information in the military documents we have. They record that he is Jewish (as opposed to "Russian", which is how basically all the other soldiers in these lists are described). Then they say that he has "higher education". One says he is an "engineer". Another is more specific, and says he's an "engineer economist" (Инж. Эконом.). They also say that he is not a member of the Communist Party.

They say he is a widower, and that his wife's name was Evdokiya Ivanovna (Евдокия Иван.). They also list his "mother", giving her name as Мария Григ. ("Maria Grig.", perhaps short for "Grigorievna"). And then they list an address: Москва С. Набер. д. 26 кв. 1ч6, which is presumably 26 Sofiyskaya Embankment, Apartment 1-6, Moscow.

Where is that address? Well, it turns out it's in the very center of Moscow ("inside the Garden Ring"), with the front looking over the Moscow River directly at the Kremlin:

Here's a current picture of the building

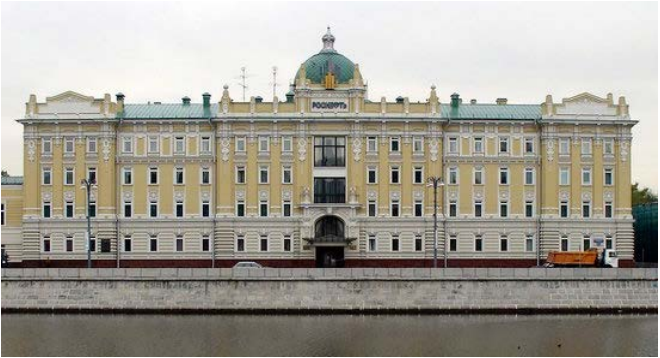

as well as one from perhaps 100 years earlier:

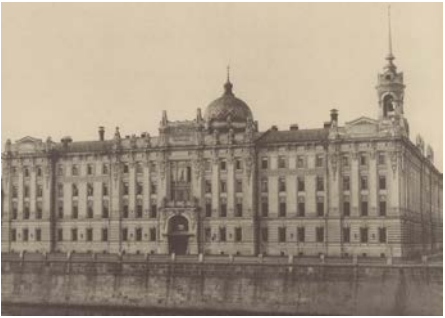

The building was built by a family of merchants named the Bakhrushins in 1900–1903 to provide free apartments for widows and orphans (apparently there were about 450 one-room 150-to-300-square-foot apartments). In the Russian Revolution, the building was taken over by the government, and set up to house the Ministry of Oil and Gas. But some "communal apartments" were left, and it's presumably in one of those that Gregory Schönfinkel lived. (Today the building is the headquarters of the Russian state oil company Rosneft.)

OK, but let's unpack this a bit further. "Communal apartments" basically means dormitory-style housing. A swank building, but apparently not so swank accommodation. Well, actually, in Soviet times dormitory-style housing was pretty typical in Moscow, so this really was a swank setup.

But then there are a couple of mysteries. First, how come a highly educated engineering economist with a swank address was just a private in the army? (When the hospital train picked up Gregory, along with four other privates, one of the others was listed as a carpenter; the others were all listed as "с/хоз" or "сельское хозяйство", basically meaning "farm laborer", or what before Soviet times would have been called "peasant").

Maybe the Russian army was so desperate for recruits after all their losses that—despite being 44 years old—Gregory was drafted. Maybe he volunteered (though then we have to explain why he didn't do that earlier). But regardless of how he wound up in the army, maybe his status as a private had to do with the fact that he wasn't a member of the Communist Party. At that time, a large fraction of the city-dwelling "elite" were members of the Communist Party (and it wouldn't have been a major problem that he was Jewish, though coming from a merchant family might have been a negative). But if he wasn't in the "elite", how come the swank address?

A first observation is that his wife's first name *Evdokiya* was a popular Russian Orthodox name, at least before 1917 (and is apparently popular again now). So presumably Gregory had—not uncommonly in the Soviet era—married someone who wasn't Jewish. But now let's look at the "mother's" name: "Мария Григ." ("Maria Grig.").

We know Gregory's (and Moses's) mother's name was Maria/"Masha" Gertsovna Schönfinkel (née Lurie)—or Мария ("Маша") Герцовна Шейнфинкель. And according to other information, she died in 1936. So—unless someone miswrote Gregory's "mother's" name— the patronymics (second names) don't match. So what's going on?

My guess is that the "mother" is actually a mother-in-law, and that it was her apartment. Perhaps her husband (most likely at that point not her) had worked at the Ministry of Oil and Gas, and that's how she ended up with the apartment. Maybe Gregory worked there too.

OK, so what was an "engineer economist" (Инженер Экономист)? In the planning-oriented Soviet system, it was something quite important: basically a person who planned and organized production and labor in some particular industry.

How did one become an "engineer economist"? At least a bit later, it was a 5-year "master's level" course of study, including courses in engineering, mathematics, bookkeeping, finance, economics of a particular sector, and "political economy" (à la Marx). And it was a very Soviet kind of thing. So the fact that that was what Gregory did presumably means that he was educated in the Soviet Union.

He must have finished high school right when the Tsar was being overthrown. Probably too late to be involved in World War I. But perhaps he got swept up in the Russian Civil War. Or maybe he was in college then, getting an early Soviet education. But, in any case, as an engineer economist it's pretty surprising that in World War II he didn't get assigned to something technical in the army, and was just a simple private in the infantry.

From the data we have, it's not clear what was going on. But maybe it had something to do with Moses.

It's claimed that Moses died in 1940 or 1942 and was "living in a communal apartment". Well, maybe that communal apartment was actually Gregory's (or at least his mother-in-law's) apartment. And here's a perhaps fanciful theory: Gregory joined the army out of some kind of despondency. His wife died. His older brother died. And in February 1942 any of his family members still in Ekaterinoslav probably died in the massacre of the Jewish population there (at least if they hadn't evacuated as a result of earlier bombing). Gregory hadn't joined the army earlier in the war, notably during the Battle of Moscow. And by 1943 he was 44 years old. So perhaps in some despondency—or anger—he volunteered for the army.

We don't know. And at this point the trail seems to go cold. It doesn't appear that Gregory had any children, and we haven't been able to find out anything more about him.

But I consider it progress that we've managed to identify that Moses's younger brother lived in Moscow, potentially providing a plausible reason that Moses might have gone to Moscow.

Actually, there may have been other "family reasons". There seems to have been quite a lot of back-and-forth in the Jewish population between Moscow and Ekaterinoslav. And Moses's mother came from the Lurie family, which was prominent not only in Ekaterinoslav, but also in Moscow. And it turns out that the Lurie family has done a fair amount of genealogy research. So we were able, for example, to reach a first cousin once removed of Moses's (i.e. someone whose parent shared a grandparent with Moses, or 1/32 of the genetics). But so far nobody has known anything about what happened to Moses, and nobody has said "Oh, and by the way, we have a suitcase full of strange papers" or anything.

I haven't given up. And I'm hoping that we'll still be able to find out more. But this is where we've got so far.

## One More Thing

In addition to pursuing the question of the fate of Moses Schönfinkel, I've made one other potential connection. Partly in compiling a bibliography of combinators, I discovered a whole collection of literature about "combinatory categorial grammars" and "combinatory linguistics".

What are these? These days, the most common way to parse an English sentence like "I am trying to track down a piece of history" is a hierarchical tree structure—analogous to the way a context-free computer language would be parsed:

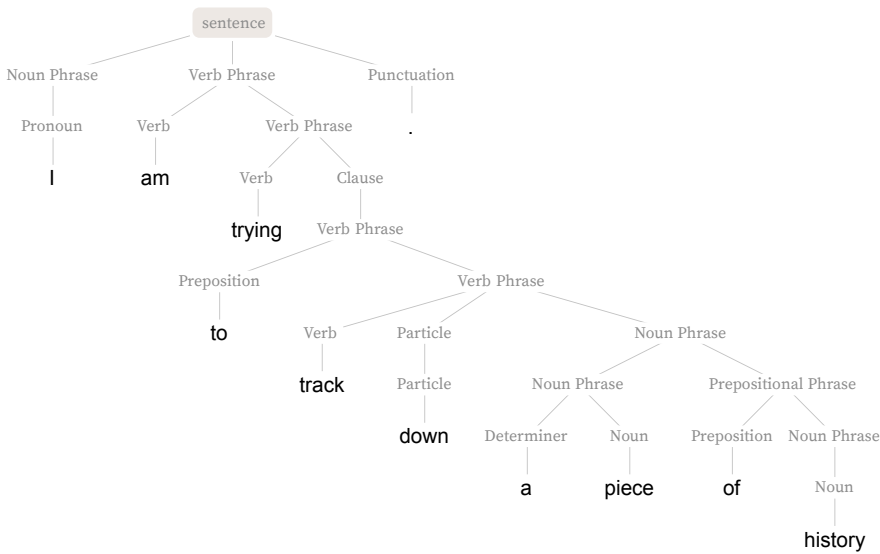

But there is an alternative—and, as it turns out, significantly older—approach: to use a so-called dependency grammar in which verbs act like functions, "depending" on a collection of arguments:

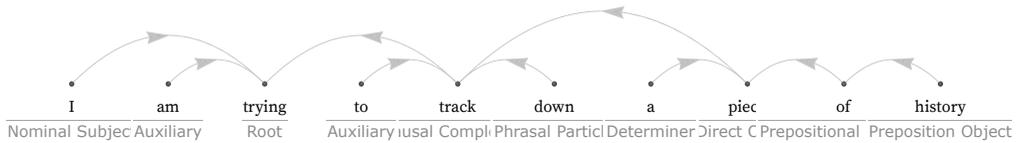

In something like Wolfram Language, the arguments in a function would appear in some definite order and structure, say as $f[x, y, z]$. But in a natural language like English, everything is just given in sequence, and a function somehow has to have a way to figure out what to grab. And the idea is that this process might work like how combinators written out in sequence "grab" certain elements to act on.

This idea seems to have a fairly tortuous history, mixed up with attempts and confusions about connecting the syntax (i.e. grammatical structure) of human languages to their semantics (i.e. meaning). The core issue has been that it's perfectly possible to have a syntactically correct sentence ("The flying chair ate a happy semicolon") that just doesn't seem to have any "real-world" meaning. How should one think about this?

I think the concept of computational language that I've spent so many years developing actually makes it fairly clear. If one can express something in computational language there's a way to compute from it. Maybe the resulting computation will align with what happens in the real world; maybe it won't. But there's some "meaningful place to go" with what one has. And the point is that a computational language has a well-defined "inner computational representation" for things. The particular syntax (e.g. sequence of characters) that one might use for input or output in the computational language is just something superficial.

But without the idea of computational language people have struggled to formalize semantics, tending to try to hang what they're doing on the detailed structure and syntax of human languages. But then what should one do about syntactically correct structures that don't "mean anything"? An example of what I consider to be a rather bizarre solution—embodied in so-called Montague grammars from the 1970s—is essentially to turn pieces of certain sentences into functions, in which there's nothing "concrete" there, just "slots" where things could go ("$x\_$ ate $y\_$")—and where one can "hold off meaninglessness" by studying things without explicitly filling in the slots.

In the original formulation, the "functions" were thought about in terms of lambdas. But combinatory categorial grammars view them instead in terms of combinators, in which in the course of a sentence words in a sense "apply to each other". And even without the notion of slots one can do "combinatory linguistics" and imagine finding the structure of sentences by taking words to "apply themselves" "across the sentence" like combinators.

If well designed (as I hope the Wolfram Language is!) computational language has a certain clean, formal structure. But human natural language is full of messiness, which has to be untangled by natural language understanding—as we've done for so many years for Wolfram|Alpha, always ultimately translating to our computational language, the Wolfram Language.

But without the notion of an underlying computational language, people tend to feel the need to search endlessly for formal structure in human natural language. And, yes, some exists. But—as we see all the time in actually doing practical natural language understanding for Wolfram|Alpha— there's a giant tail that seems to utterly explode any all-encompassing formal theory.

Are there at least fragments that have formal structure? There are things like logic ("and", "or", etc.) that get used in human language, and which are fairly straightforwardly formalizable. But maybe there are more "functional" structures too, perhaps having to do with the operation of verbs. And in combinatory linguistics, there've been attempts to find these—even for example directly using things like Schönfinkel's $S$ combinator. (Given $S f g x \rightarrow f[x][g[x]]$ one can start imagining—with a slight stretch—that "eat peel orange" operates like the $S$ combinator in meaning "eat[orange][peel[orange]]".)

Much of the work on this has been done in the last few decades. But it turns out that its history stretches back much further, and might conceivably actually intersect with Moses Schönfinkel himself.

The key potential link is Kazimierz Ajdukiewicz (1890–1963). Ajdukiewicz was a Polish logician/philosopher who long tried to develop a "mathematicized theory" of how meaning emerges, among other things, from natural language, and who basically laid the early groundwork for what's now combinatory linguistics.

Kazimierz Ajdukiewicz was born two years after Moses Schönfinkel, and studied philosophy, mathematics and physics at the University of Lviv (now in Ukraine), finishing his PhD in 1912 with a thesis on Kant's philosophy of space. But what's most interesting for our purposes is that in 1913 Ajdukiewicz went to Göttingen to study with David Hilbert and Edmund Husserl.

In 1914 Ajdukiewicz published one paper on "Hilbert's New Axiom System for Arithmetic", and another on contradiction in the light of Bertrand Russell's work. And then in 1915 Ajdukiewicz was drafted into the Austrian army, where he remained until 1920, after which he went to work at the University of Warsaw.

But in 1914 there's an interesting potential intersection. Because June of that year is when Moses Schönfinkel arrived in Göttingen to work with Hilbert. At the time, Hilbert was mostly lecturing about physics (though he also did some lectures about "principles of mathematics"). And it seems inconceivable that—given their similar interests in the structural foundations of mathematics—they wouldn't have interacted.

Of course, we don't know how close to combinators Schönfinkel was in 1914; after all, his lecture introducing them was six years later. But it's interesting to at least imagine some interaction with Ajdukiewicz. Ajdukiewicz's own work was at first most concerned with things like the relationship of mathematical formalism and meaning. (Do mathematical constructs "actually exist", given that their axioms can be changed, etc.?) But by the beginning of the 1930s he was solidly concerned with natural language, and was soon writing papers with titles like "Syntactic Connexion" that gave formal symbolic descriptions of language (complete with "functors", etc.) quite reminiscent of Schönfinkel's work.

So far as I can tell Ajdukiewicz never explicitly mentioned Schönfinkel in his publications. But it seems like too much of a coincidence for the idea of something like combinators to have arisen completely independently in two people who presumably knew each other—and never to have independently arisen anywhere else.

### Thanks

Thanks to Vitaliy Kaurov for finding additional documents (and to the State Archives of the Dnipropetrovsk Region and Elena Zavoiskaia for providing various documents), Oleg and Anna Marichev for interpreting documents, and Jason Cawley for information about military history. Thanks also to Oleg Kiselyov for some additional suggestions on the original version of this piece.

# References

*Links to references are included within the body of this document.*